\documentclass[11pt]{article}

\usepackage[margin=1in]{geometry}
\usepackage{parskip}

\usepackage{amsmath}
\usepackage{amssymb}

\usepackage{graphicx}
\usepackage{booktabs}
\usepackage{needspace}
\usepackage{float}
\usepackage{algorithm}
\usepackage{algpseudocode}
\usepackage{cite}

\usepackage{placeins}

\usepackage{microtype}
\usepackage{hyperref}
\usepackage{xcolor}

\title{
In-Context Learning for Multi-Fidelity Multiscale Model Closure and Cross-Operator Out-Of-Distribution Generalization
}

\author{
Xinyuan Li\thanks{Department of Applied and Computational Mathematics
and Statistics, University of Notre Dame, Notre Dame, IN 46556, USA.
E-mail: \href{mailto:xli46@nd.edu}{xli46@nd.edu}.}
\and
Matthew Zahr\thanks{Department of Aerospace and Mechanical Engineering,
University of Notre Dame, Notre Dame, IN 46556, United States.
E-mail: \href{mailto:mzahr@nd.edu}{mzahr@nd.edu}.}
\and
Panos Stinis\thanks{Advanced Computing, Mathematics, and Data Division,
Pacific Northwest National Laboratory, Richland, WA 99354, USA.
E-mail: \href{mailto:Panagiotis.stinis@pnnl.gov}{Panagiotis.stinis@pnnl.gov}.}
\and
Zecheng Zhang\thanks{Corresponding author. Department of Applied and
Computational Mathematics and Statistics, University of Notre Dame,
Notre Dame, IN 46556, USA.
E-mail: \href{mailto:zzhang48@nd.edu}{zzhang48@nd.edu} or \href{mailto:zzhang48@nd.edu}{zecheng.zhang.math@gmail.com}.}
}

\date{}

\begin{document}

\maketitle


\begin{abstract}


Multi-fidelity learning uses low-fidelity information (coarse-scale
PDE models) together with limited high-fidelity data (fine-scale
solutions) to obtain high-accuracy PDE solutions. It is a data-driven
approach to downscale reduced-order models to solve PDEs efficiently
and accurately.
Standard multi-fidelity learning focuses on one operator or PDE,
and, given a new PDE, retraining the model using a large set of data
(coarse--fine solution pairs of the new PDE) may be required,
which is costly.
In this work, we propose an in-context multi-fidelity learning
framework for learning multiple multi-fidelity operators with one
shared model.
In-context learning aims to learn to downscale a coarse-scale
model to a fine-scale model using a few example coarse--fine pairs.
Compared to the standard model, the proposed method is able to
construct fine-scale predictions for unseen operators with only
a few examples (at low cost), or even no additional examples if
the testing operators match the training operators but have a
different prediction goal, without fine-tuning.
Specifically, the model is trained using a small number of paired
low-fidelity (coarse PDE solution) and high-fidelity (fine PDE solution)
examples as demonstrations from multiple operators, together with
a new coarse-scale solution query, to predict the fine-scale solution
corresponding to the query.
Once trained, given a new PDE operator (e.g., an initial-condition-to-solution
map), a few (fewer than five) coarse--fine solution pairs as
demonstrations, and a target query consisting of the low-cost coarse
solution for the query initial condition, the model can predict
the corresponding fine-scale solution without any model fine-tuning.
Notably, during testing, if the target operator appears in the
training dataset but has a different prediction goal, e.g.,
predicting the solution beyond the training time interval, our model
does not need to generate more demonstrations. Instead, the model
can leverage the training coarse--fine pairs from this PDE and
information from other training operators to make the prediction
without fine-tuning. The model therefore has a very low cost during
testing and only requires generating a low-cost coarse-scale
solution query.
Notably, our smallest model contains only 255,364 trainable
parameters, corresponding to approximately 456 parameters per
training operator for 560 operators.

\end{abstract}


\section{Introduction}

Operator learning \cite{lu2021learning, li2020fourier,lin2023b,zhang2023d2no,lin2023learning,wen2022u,hao2025laplacian,yu2024nonlocal,chiu2025context,sanderse2024scientific,meuris2023machine, xu2023transfer,molinaro2023neuralinverseoperatorssolving,Nelsen_2026,zhang2026numerical,Long2024Kernel} aims to approximate mappings between function spaces and has become an important approach for learning solution operators of differential equations. Early studies established universal approximation results for continuous functionals and nonlinear operators using neural networks \cite{chen1993approximations,chen1995universal, back2002universal}. 
Subsequent work has provided further theoretical understanding of neural operators, including 
the convergence of the operator approximation error with respect to network size \cite{liu2024neural,yang2026generalization,bayraktar2025deep, liu2024deep,lanthaler2023parametric,lanthaler2022error}.
Operator learning has also been extended to multi-operator learning \cite{back2002universal}, which goes beyond the setting of learning a single operator. 
Extending single-operator approaches to multiple operators can become costly when separate models or repeated training are required for different operators. 
To address this limitation, recent multi-operator learning approaches \cite{zhang2024modno,sun2025towards, cao2024vicon, yang2023context, liu2023prose, ovadia2023vitovisiontransformeroperator, yu2024nonlocal} use shared neural architectures to learn multiple operators within a single model, allowing information to be shared across related operator-learning tasks \cite{zhang2024modno,sun2025towards}. 
A key advantage of such approaches is the ability to generalize to previously unseen operators at test time when these operators are drawn from a distribution similar to that of the training operators.

Multi-fidelity learning \cite{lu2022multifidelity, howard2023multifidelity, zhang2024bayesian,sanderse2024scientific, xu2026improving,Ahmed_2023} combines models with different levels of accuracy and computational cost. 
Low-fidelity solutions are typically cheaper to compute, while high-fidelity solutions are more accurate but more expensive. 
The goal is to use low-fidelity information together with limited high-fidelity data to improve prediction accuracy while reducing computational cost.
A variety of approaches have been developed to model the relation between low- and high-fidelity data, including statistical, nonlinear, and neural-network-based methods.

Unlike model reduction \cite{Zahr2015Progressive, Amsallem2012Nonlinear, Amsallem2015Fast, Washabaugh2012Nonlinear} or multiscale methods \cite{efendiev2000convergence, Efendiev_2013, galvis2010domain, Schaeffer2020Extracting, Mackey2014Compressive, Hardwick2025}, which compute a coarse-scale approximation to the fine-scale model, multi-fidelity learning downscales the coarse-scale solution to the fine scale using data.
Hence, multi-fidelity learning provides a data-driven approach to computing fine-scale model solutions by improving the accuracy of solutions obtained from low-cost coarse-scale numerical simulations.
However, data are not easy to obtain, so the goal is to train a model that can be reused for new PDEs or for the same PDE under different settings, such as times beyond the training range.

Multi-fidelity learning has also been combined with neural operators. Existing approaches include residual or discrepancy learning \cite{lu2022multifidelity,de2023bi}, composite multi-fidelity operator models \cite{howard2023multifidelity}, and direct coarse-to-fine (multiscale downscaling) mappings \cite{zhang2024bayesian}. These methods show that low-fidelity solutions can be used effectively to reduce the need for expensive high-fidelity data.
However, many multi-fidelity operator-learning methods are developed for a specific coarse-to-fine relation or a single operator-learning task, making extension to multiple operators costly.
Beyond this limitation, out-of-distribution (OOD) prediction presents another challenge, in which the testing query lies outside the input range covered by the training data for the target PDE. The coarse-to-fine relation learned within the training range may not remain accurate in a new setting, such as later times or different initial conditions.

In-context learning (ICL) \cite{brown2020language,garg2022can,xie2021explanation,von2023transformers,akyurek2022learning, cole2026context, cole2026context1, cole2024context, yang2023context, yang2024pde} uses a small number of input--output examples as demonstrations to help the model identify the relation relevant to a new query while keeping the model parameters fixed.
Popularized by large language models \cite{brown2020language}, ICL has also been studied for learning functions from input--output examples \cite{garg2022can}. More recently, this idea has been extended to operator learning. ICON \cite{yang2023context} uses prompted solution pairs to identify the target operator during inference, allowing one shared model to handle multiple operators without retraining or fine-tuning. 
Subsequent work has studied generalization to unseen PDE operators \cite{yang2024pde}. 
VICON \cite{cao2024vicon} extends in-context operator learning to dense spatial data in multi-physics problems using a vision-based architecture and also considers flexible temporal mappings with varying timestep strides. 
FMint \cite{song2025fmint} applies in-context learning to numerical correction (differences) for ODEs. It uses a large model with approximately 15.8 million trainable parameters, and its reported OOD evaluations involve additional fine-tuning on trajectories from the target system.
The ability to use demonstrations to identify the relation relevant to a new query while keeping the model parameters fixed makes ICL well suited to multi-fidelity learning. 
In this setting, the shared model learns a family of coarse-to-fine relations, while the demonstrations provide the context needed to identify the relation relevant to the target operator. Once trained, the same model can adapt to different operators through the supplied demonstrations without parameter updates. This provides a natural way to combine multi-fidelity learning with multi-operator generalization. Although ICL is often associated with large-scale models, the multi-fidelity setting can support a more compact model because the coarse and fine solutions within each demonstration share similar structure.

In this work, we propose an in-context multi-fidelity learning framework for learning a family of coarse-to-fine operators with one shared model. Once trained, the model parameters remain fixed, and no fine-tuning is required at inference. The framework generalizes to unseen operators drawn from the training distribution and also supports operator-parameter OOD extrapolation beyond the training range. It further handles OOD query regimes, including initial-condition and temporal ranges not represented for the target operators during training. At inference, the cost is already low because only a few demonstrations are needed and no parameter updates are required. 
More importantly, we further reduce the inference cost for OOD tasks through cross-operator transfer, where OOD information learned from a group of operators with more training data improves prediction for another group with limited training data that has not observed the same OOD regime. 
In this setting, demonstrations used during training are reused for the new OOD queries. As a result, a new query requires only the coarse solution generated by the low-fidelity solver, with no new fine solution or fine-tuning. This is particularly advantageous in multi-fidelity settings, where low-fidelity simulations are inexpensive while high-fidelity simulations are costly.

The framework is also computationally compact. Our smallest model contains 255,364 trainable parameters, corresponding to approximately 456 parameters per training operator for 560 operators, and the other models used in this work are also relatively small-scale models. The model can accommodate different numbers of demonstrations without changing its architecture or parameter count, and we study the effect of both the number and consistency of the supplied demonstrations. Beyond standard same-time coarse-to-fine mappings, we also consider temporally shifted mappings of the form \(v_c(t)\mapsto v_f(t+\Delta t)\). We further extend the framework to a two-dimensional multi-fidelity problem and compare it with existing multi-fidelity operator-learning methods. Together, these experiments show that the shared-model framework can adapt across unseen operators and multiple OOD regimes through in-context demonstrations, including cross-operator transfer to regimes not observed for the target operators during training.


\section{Methodology}

We consider an in-context multi-fidelity learning framework for a
family of operators within a given PDE application. Rather than
training a separate model for each operator, we use a single model with
shared trainable parameters across the operator family.

\begin{figure}[H]
    \centering
    \includegraphics[width=\textwidth]{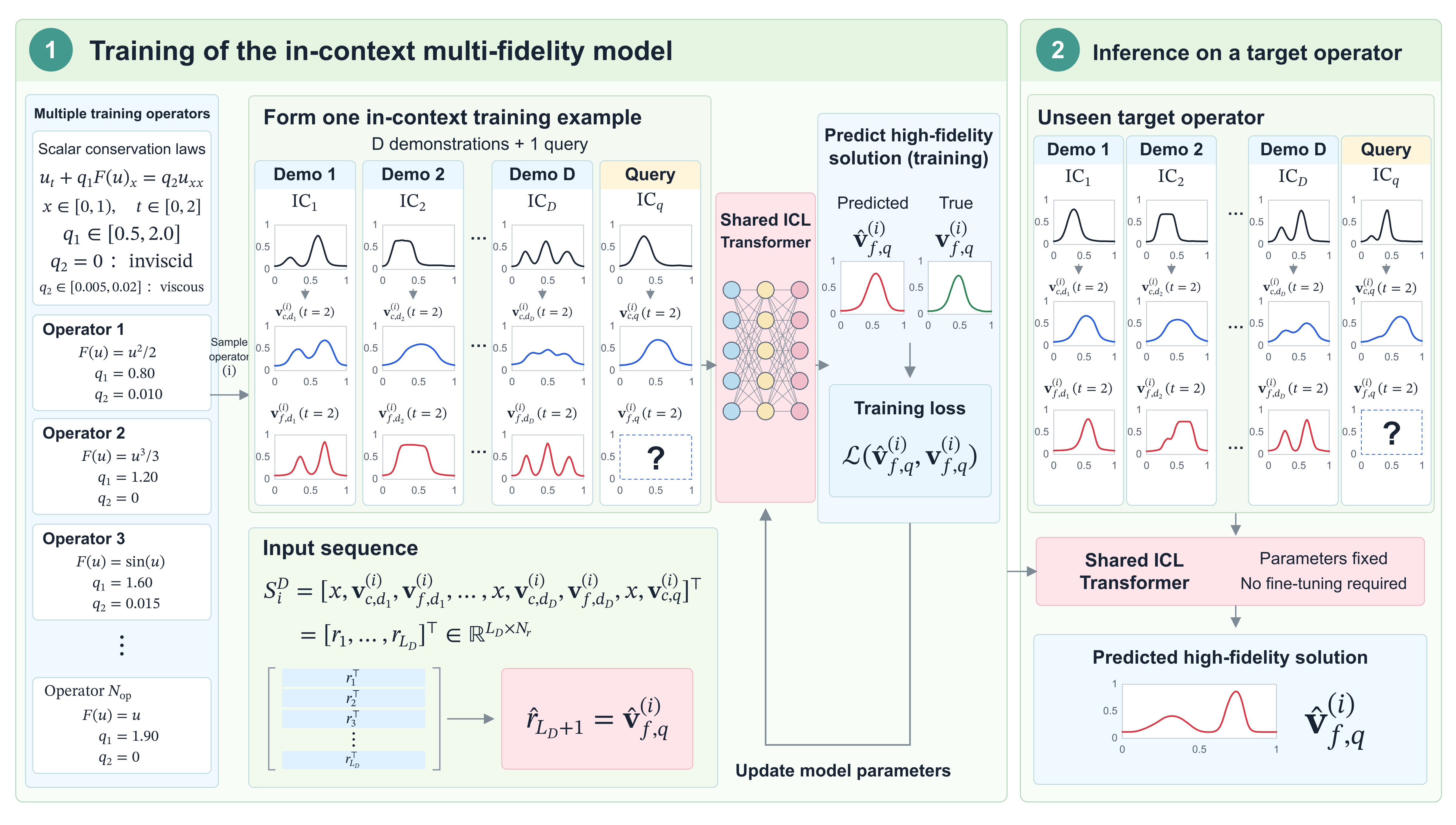}
    \caption{
    Overview of the proposed in-context multi-fidelity learning
    framework. Paired coarse--fine demonstrations from a target
    operator and a coarse query are arranged into an ordered sequence.
    A Transformer predicts the missing fine-solution row corresponding
    to the query through a one-row shift, while the same shared
    trainable parameters are used across the operator family.
    }
    \label{fig:method_overview}
\end{figure}

For each prediction, a small number of paired coarse--fine
demonstrations from the target operator are provided together with a
new coarse query. These demonstrations provide the context for the
multi-fidelity operator associated with the target operator. We arrange
the demonstrations and query into an ordered sequence and treat the
task as a next-row prediction problem. The sequence is processed by a
Transformer, and the fine solution corresponding to the query is
predicted as the missing next row. The number of demonstrations can
vary, while the shared trainable parameters remain fixed during
inference. An overview of the framework is shown in
Figure~\ref{fig:method_overview}.

\subsection{Multi-fidelity operator formulation}
\label{sec:mf_operator_formulation}

Let $U_i$ denote the $i$-th input function space with domain $\Omega_{u,i}\subset\mathbb R^{d_{u,i}}$ and let $V_i$ denote the
corresponding space of output functions defined on
$\Omega_{v,i}\subset\mathbb R^{d_{v,i}}$.
Here, $d_{u,i}$ and $d_{v,i}$ denote the dimensions of the
input and output domains, respectively.
We consider $N_{\mathrm{op}}$ operators
$G_i:U_i\rightarrow V_i$, $i=1,\ldots,N_{\mathrm{op}}$.
For a given input $u\in U_i$, the corresponding output function is
$v^{(i)}=G_i(u)\in V_i$. 

In the multi-fidelity setting, $G_i$ is the \textit{exact} solution operator of the PDEs. For example, it maps the initial condition $u\in U_i$ to the solution of the PDEs $v^{(i)}=G_i(u)$.
The solution operator can be evaluated using numerical solvers
at low- and high-fidelity levels. 
The two levels differ in numerical
accuracy and computational cost, for example through differences in
discretization resolution, orders of accuracy, or other solver settings.
At the low-fidelity level, we denote the \textit{coarse-scale} \textit{numerical} solution operator as
$G_i^c:U_i\rightarrow V_i^c$, which produces the coarse solution
$v_c^{(i)} = G_i^c(u)$, where $V_i^c$ denotes the corresponding
coarse-solution output space.
At the high-fidelity level, we denote the \textit{fine-scale} \textit{numerical}
solution operator as $G_i^f:U_i\rightarrow V_i^f$, which produces the fine
solution $v_f^{(i)} = G_i^f(u)$, where $V_i^f$ denotes the
corresponding fine-solution output space.

Let $\mathcal A_i$ denote the space of application-specific
information associated with inputs in $U_i$. For each input
$u\in U_i$, the map $a:U_i\rightarrow\mathcal A_i$ provides
the auxiliary information $a(u)$ used in the coarse-to-fine prediction.
We define the corresponding
\textit{multi-fidelity} operator
$G_i^{c\rightarrow f}:\mathcal A_i\times V_i^c\rightarrow V_i^f$
and write
$v_f^{(i)}=G_i^{c\rightarrow f}\bigl(a(u),v_c^{(i)}\bigr)$
for $i=1,\ldots,N_{\mathrm{op}}$.
The role and sequence representation of $a(u)$ are described
in Section~\ref{sec:icl_framework}, and the application-specific forms
are summarized in Appendix~\ref{app:sequence_summary}.
In this work, we consider the family of multi-fidelity operators
$\{G_i^{c\rightarrow f}\}_{i=1}^{N_{\mathrm{op}}}$.

Standard multi-fidelity operator learning considers a single
multi-fidelity operator $G_i^{c\rightarrow f}$ over a specified
training input range
\cite{howard2023multifidelity,zhang2024bayesian} for one PDE solution operator $G_i$.
For a fixed PDE problem, the trained model predicts fine solutions
from coarse solutions. 
When the PDE problem changes (different $i$), additional
re-training or fine-tuning with additional data is needed to learn the new multi-fidelity
operator, limiting its use across different PDE problems. 
In this work, we learn multiple multi-fidelity operators with one
shared model and address OOD prediction without test-time fine-tuning.

To describe the multi-operator learning and OOD prediction settings,
let $\mathcal I_{\mathrm{train}}$ and $\mathcal I_{\mathrm{test}}$
denote the index sets of the multi-fidelity operators
$G_i^{c\rightarrow f}$ used for training and held out for testing,
respectively, with
$\mathcal I_{\mathrm{train}},
\mathcal I_{\mathrm{test}}
\subset\{1,\ldots,N_{\mathrm{op}}\}$ and
$\mathcal I_{\mathrm{train}}\cap
\mathcal I_{\mathrm{test}}=\varnothing$.
For each $i\in\mathcal I_{\mathrm{train}}$, let
$\mathcal U_i^{\mathrm{train}}\subset U_i$ denote the input range
from which training inputs are sampled.

We consider three different extrapolation settings, including
prediction outside the training parameter range.
Notably, standard single-operator multi-fidelity learning may
require costly retraining using additional data.
All three settings are tested in Section~\ref{sec:numerical_examples}.
(1) For unseen-operator prediction, the testing multi-fidelity operator
$G_i^{c\rightarrow f}$ has index $i\in\mathcal I_{\mathrm{test}}$
and is not used during training.
(2) For query OOD prediction, the
target multi-fidelity operator $G_i^{c\rightarrow f}$ has index
$i\in\mathcal I_{\mathrm{train}}$, but the query input lies outside
$\mathcal U_i^{\mathrm{train}}$.
(3) For operator-parameter OOD prediction,
the target multi-fidelity operator $G_i^{c\rightarrow f}$ is held
out during training, and the parameters defining the corresponding
solution operator $G_i$ lie outside the training parameter ranges.


\subsection{In-context learning framework}
\label{sec:icl_framework}

We formulate an in-context learning framework for the family of
multi-fidelity operators
$\{G_i^{c\rightarrow f}\}_{i=1}^{N_{\mathrm{op}}}$,
where $i$ indexes the \textit{exact} PDE solution operator $G_i$ and its corresponding
multi-fidelity operator $G_i^{c\rightarrow f}$.
The goal is to learn multiple multi-fidelity operators using one
shared model.
A paired coarse--fine demonstration
$(v_c^{(i)} = G_i^c(u), v_f^{(i)} = G_i^f(u))$ consists of the coarse and fine
solutions generated from the same input $u$.
These demonstrations provide the context needed to identify the
target multi-fidelity operator. Once trained, the shared model can
be applied across different target multi-fidelity operators while
its parameters remain fixed, with no retraining or fine-tuning
at inference.
We describe the training and testing procedures separately. See Figure \ref{fig:method_overview} for the demonstration.

\paragraph{Training.}

For a target multi-fidelity operator $G_i^{c\rightarrow f}$,
let $D$ denote the number of demonstrations and let
$d_1,\ldots,d_D$ be their distinct sample indices in sequence order.
During training, the target operator has index
$i\in\mathcal I_{\mathrm{train}}$, and the demonstration and query
inputs are selected from $\mathcal U_i^{\mathrm{train}}$.
The notation and sequence construction in \eqref{eq:icl_sequence} also apply at inference.

For $j=1,\ldots,D$, the demonstration input
$u_{d_j}^{(i)}\in U_i$ generates the coarse--fine pair
$(v_{c,d_j}^{(i)},v_{f,d_j}^{(i)})$,
where
$v_{c,d_j}^{(i)}=G_i^c(u_{d_j}^{(i)})$ and
$v_{f,d_j}^{(i)}=G_i^f(u_{d_j}^{(i)})$.
For example, the input may be an initial condition for a
terminal-time solution or a time associated with a solution trajectory.
Evaluating the coarse and fine solution operators at this input
provides the paired solutions used as a demonstration.

For a new query with sample index
$q\notin\{d_1,\ldots,d_D\}$ and input
$u_q^{(i)}\in U_i$, we denote the corresponding solutions as
$v_{c,q}^{(i)}=G_i^c(u_q^{(i)})$ and
$v_{f,q}^{(i)}=G_i^f(u_q^{(i)})$.
For example, the query input may be a different initial condition
or another time along the trajectory, evaluated using the same
solution operators as the demonstrations.
The fine query solution $v_{f,q}^{(i)}$ is withheld from the model
input and is used as the supervision target during training.

Before constructing the model input, we define the discretization of the coarse and fine solutions. To simplify notation, we omit the operator and input indices in this paragraph. Let
$\widetilde{\mathbf v}_c\in\mathbb R^{N_c}$ and
$\widetilde{\mathbf v}_f\in\mathbb R^{N_f}$
denote the coarse and fine solutions on the original coarse and fine solver grids, with $n_c$ and $n_f$ spatial locations, respectively, and $n_v$ physical components.
$N_c=n_vn_c$ and $N_f=n_vn_f$ denote the corresponding numbers of entries.
For a scalar-valued coarse solution,
$\widetilde{\mathbf v}_c=[v_c(x_1),\ldots,v_c(x_{n_c})]^\top$
and $N_c=n_c$.
Within each numerical application, let $n_r$ denote the number of spatial locations on the common learning grid. We interpolate
$\widetilde{\mathbf v}_c$ and $\widetilde{\mathbf v}_f$
from the respective solver grids onto the common learning grid, obtaining
$\mathbf v_c,\mathbf v_f\in\mathbb R^{N_r}$,
where $N_r=n_vn_r$ denotes the number of entries in one solution row on the common learning grid.
Each solution row contains the values of all $n_v$ physical components at the $n_r$ spatial locations.
The coarse and fine rows use the same ordering of spatial locations and physical components.
We use bold symbols for the corresponding finite-dimensional solution rows.
For $j=1,\ldots,D$,
$\mathbf v_{c,d_j}^{(i)},\mathbf v_{f,d_j}^{(i)}\in\mathbb R^{N_r}$,
and
$\mathbf v_{c,q}^{(i)},\mathbf v_{f,q}^{(i)}\in\mathbb R^{N_r}$.

The application-specific information $a(u)$ is included when needed
to improve the conditioning of the coarse-to-fine prediction.
For a sample with index $s\in\{d_1,\ldots,d_D,q\}$, let
$a_s^{(i)}=a(u_s^{(i)})$ denote this information, so that
$v_{f,s}^{(i)}
=
G_i^{c\rightarrow f}
\bigl(a_s^{(i)},v_{c,s}^{(i)}\bigr)$.
Suppose $a_s^{(i)}$ is represented by $n_a$ auxiliary sequence
rows for each sample. We denote the finite-dimensional auxiliary-row
representations by bold symbols. After any required interpolation or
broadcasting to the common row length,
$\mathbf a_{s,m}^{(i)}\in\mathbb R^{N_r}$ for
$m=1,\ldots,n_a$.
Each auxiliary row and each coarse or fine solution row corresponds
to one token. We use $\boldsymbol{\alpha}_s^{(i)}$ as shorthand
for the ordered list of auxiliary rows
$\mathbf a_{s,1}^{(i)},\ldots,\mathbf a_{s,n_a}^{(i)}$.
In the sequence in \eqref{eq:icl_sequence}, each $\boldsymbol{\alpha}_s^{(i)}$
expands into these $n_a$ row vectors before concatenation.
The application-specific auxiliary-row quantities are summarized
in Table~\ref{tab:app_sequence_summary} in Appendix~\ref{app:sequence_summary}.

The $D$ demonstrations and the query are arranged into a single
ordered sequence $S_i^D$,
\begin{equation}
S_i^D=
[\boldsymbol{\alpha}_{d_1}^{(i)},
\mathbf v_{c,d_1}^{(i)},\mathbf v_{f,d_1}^{(i)},
\ldots,
\boldsymbol{\alpha}_{d_D}^{(i)},
\mathbf v_{c,d_D}^{(i)},\mathbf v_{f,d_D}^{(i)},
\boldsymbol{\alpha}_q^{(i)},\mathbf v_{c,q}^{(i)}]^\top
=[r_1,\ldots,r_{L_D}]^\top
\in\mathbb R^{L_D\times N_r},
\label{eq:icl_sequence}
\end{equation}
where $r_\ell\in\mathbb R^{N_r}$ is the vector whose transpose
forms the $\ell$th row of the sequence, and
$L_D=D(n_a+2)+(n_a+1)$ is the total number of valid input rows.
Each demonstration contributes $n_a+2$ rows, and the query contributes
$n_a+1$ rows because the fine query solution is withheld.
Different selections of demonstrations and query produce different
input sequences even when $i$ and $D$ are unchanged.

For an example with $D$ demonstrations, the final valid input row is
the coarse query,
$r_{L_D}=\mathbf v_{c,q}^{(i)}$, and the fine query target is treated
as the missing next row,
$r_{L_D+1}=\mathbf v_{f,q}^{(i)}$. The prediction follows the one-row
shift
\begin{equation}
\begin{bmatrix}
r_1^\top\\
r_2^\top\\
\vdots\\
r_{L_D}^\top
\end{bmatrix}
\longrightarrow
\widehat r_{L_D+1},
\quad
\widehat r_{L_D+1}
=
\widehat{\mathbf v}_{f,q}^{(i)}.
\label{eq:icl_one_row_shift}
\end{equation}
The model predicts the next row $\widehat r_{L_D+1}$,
corresponding to the fine query solution
$\widehat{\mathbf v}_{f,q}^{(i)}$.

The prediction in \eqref{eq:icl_one_row_shift} is a single
autoregressive-style next-row prediction.
$\widehat{\mathbf v}_{f,q}^{(i)}$ is not fed back into the sequence for
any further prediction, so no recursive error propagation is
introduced.

Let $\mathcal M_\theta$ denote the ICL model with trainable
parameters $\theta$.
We write the ICL prediction as
\begin{equation*}
\widehat{\mathbf v}_{f,q}^{(i)}
=
\mathcal M_\theta(S_i^D)
\approx
\mathbf v_{f,q}^{(i)}.
\end{equation*}
The trainable parameters $\theta$ are shared across the family
$\{G_i^{c\rightarrow f}\}_{i=1}^{N_{\mathrm{op}}}$ and remain fixed
when the trained model is applied to a target multi-fidelity operator
$G_i^{c\rightarrow f}$ at inference.
The demonstrations in $S_i^D$ provide context for the target
multi-fidelity operator without updating $\theta$.

The number of demonstrations $D$ need not be fixed, so the input
sequence $S_i^D$ can have different lengths $L_D$. 
This motivates the use
of a Transformer, which is a natural choice for handling variable-length
input sequences without changing the number of trainable model parameters.
Further details are provided in Appendix~\ref{app:transformer_details}.

Let $\mathcal D$ denote the set of demonstration counts used
by the model. For each training example, $D$ is sampled independently
from $\mathcal D$, and the input sequence $S_i^D$ is paired with
the fine query target $\mathbf v_{f,q}^{(i)}$.

We denote the training objective by $\mathcal L(\theta)$. Its
definition is provided in Appendix~\ref{app:training_objective}.
The training procedure is summarized in Algorithm~\ref{alg:icl_training}.

\begin{algorithm}[H]
\caption{Training of the in-context multi-fidelity model.
See Figure \ref{fig:method_overview} for the demonstration.}
\label{alg:icl_training}
\begin{algorithmic}[1]

\Statex \textbf{Input:}
Paired training data, training operator-index set
$\mathcal I_{\mathrm{train}}$,
demonstration-count set $\mathcal D$,
number of updates $N_{\mathrm{iter}}$, and shared model
$\mathcal M_\theta$

\Statex \textbf{Output:}
Trained shared parameters $\theta$

\State Initialize $\theta$

\For{$n=1,\ldots,N_{\mathrm{iter}}$}
    \State Sample $i\in\mathcal I_{\mathrm{train}}$
    and $D\in\mathcal D$
    \State Select $D$ distinct demonstrations and one additional
    query from the training data for
    $G_i^{c\rightarrow f}$
    \State Construct $S_i^D$ according to
    \eqref{eq:icl_sequence}, using the demonstrations,
    auxiliary information, and coarse query
    \State Predict the fine query
    $\widehat{\mathbf v}_{f,q}^{(i)}
    =\mathcal M_\theta(S_i^D)$
    from the final valid input position
    \State Compute the loss using
    $\widehat{\mathbf v}_{f,q}^{(i)}$
    and the true fine query $\mathbf v_{f,q}^{(i)}$
    \State Update $\theta$ using the optimizer
\EndFor

\State \Return $\theta$

\end{algorithmic}
\end{algorithm}

\paragraph{Testing.}

During testing, we consider unseen-operator generalization and
cross-operator OOD transfer.
For unseen-operator generalization, the target multi-fidelity operator
$G_i^{c\rightarrow f}$ has index $i\in\mathcal I_{\mathrm{test}}$,
and few-shot demonstrations are provided for this held-out operator.
Notably, an OOD scenario can still arise even when the testing operator has been seen during training.
For example, $G_i$ maps initial conditions to solutions at times $t<T$ during training, while testing considers solutions at times $t>T$, so that the new query satisfies
$u_q^{(i)}\in U_i\setminus\mathcal U_i^{\mathrm{train}}$.
In this case, few-shot demonstrations (initial conditions and the corresponding solutions in the OOD time range) may not be provided.
Since learning can be improved by using information from other operators $G_k$, with $k\in\mathcal I_{\mathrm{train}}$ and $k\neq i$, whose training data cover $t>T$, we refer to this setting as cross-operator OOD transfer.
In this case, the demonstrations associated with the target multi-fidelity operator $G_i^{c\rightarrow f}$ are reused from training.
No new fine query solution $\mathbf v_{f,q}^{(i)}$ or fine-tuning is required.
Only the new coarse query $\mathbf v_{c,q}^{(i)}$ and the application-specific query rows $\mathbf a_{q,1}^{(i)},\ldots,\mathbf a_{q,n_a}^{(i)}$ are needed in addition to the reused demonstrations.

At inference, demonstrations and a new coarse query for the target
multi-fidelity operator $G_i^{c\rightarrow f}$ are arranged into
$S_i^D$. The trained model predicts the next row,
$\widehat{\mathbf v}_{f,q}^{(i)}
=
\mathcal M_\theta(S_i^D)$,
while the shared parameters $\theta$ remain fixed.
The data-generation cost at inference for \eqref{eq:icl_sequence}
consists of computing a few coarse--fine pairs
$(\mathbf v_{c,d_j}^{(i)},\mathbf v_{f,d_j}^{(i)})$
as few-shot demonstrations and the new coarse query
$\mathbf v_{c,q}^{(i)}$, which is inexpensive to obtain
using a low-cost numerical solver.
Notably, new few-shot demonstrations are not needed in the cross-operator OOD transfer setting described in Section~\ref{sec:icl_framework}, since the demonstrations for the target operator are reused from training.
Additional implementation details, including the
model architecture, variable-length batching, training loss, data
normalization, and application-specific training settings, are
provided in Appendices~\ref{app:common_details}--\ref{app:ice_details}.

\section{Numerical examples}
\label{sec:numerical_examples}

In this section, we apply our in-context multi-fidelity framework to four PDE applications to evaluate the generalization ability of the proposed framework across different PDE settings. The main question is whether, given only a few paired coarse--fine demonstrations, the model can predict the corresponding fine solution for a new coarse query without any test-time parameter updates. In particular, we are interested in cross-operator generalization, out-of-distribution (OOD) transfer, and the role of the demonstration context.

We use the notation defined in Sections~\ref{sec:mf_operator_formulation} and~\ref{sec:icl_framework}. 
For readability, the operator
index $i$ is omitted from solution rows in this section unless needed.
The coordinate and other application-specific inputs are represented
by the auxiliary rows defined in Section~\ref{sec:icl_framework}.

Prediction accuracy is evaluated using the relative $L^2$ error $\|\widehat{\mathbf v}_f-\mathbf v_f\|_2/\|\mathbf v_f\|_2$ and the coefficient of determination $R^2=1-\sum_{s=1}^{N_r}(\widehat v_{f,s}-v_{f,s})^2/\sum_{s=1}^{N_r}(v_{f,s}-\bar v_f)^2$, where $\bar v_f=(1/N_r)\sum_{s=1}^{N_r}v_{f,s}$. For the ice-sheet benchmark, we also report the mean-squared error
$\mathrm{MSE}=(1/N_r)\sum_{s=1}^{N_r}(\widehat v_{f,s}-v_{f,s})^2$,
where the sum includes all spatial locations and physical components.
All metrics are computed separately for each test case.
The reported MSE is averaged over the test cases and is computed
in physical velocity units, giving units of $(\mathrm{m}/\mathrm{yr})^2$. The low-fidelity baseline error is computed by replacing $\widehat{\mathbf v}_f$ with the corresponding coarse solution $\mathbf v_c$. Unless otherwise stated, results are reported as mean $\pm$ standard deviation over the complete set of test cases.

Throughout this section, one test case refers to one prediction. The model receives $D$ coarse--fine demonstrations and one coarse query $\mathbf v_{c,q}$ from the same target operator, and predicts the corresponding fine solution $\widehat{\mathbf v}_{f,q}$. The true fine solution $\mathbf v_{f,q}$ is used only to evaluate the prediction error. Thus, one operator can provide multiple test cases through different queries.

When the target operator is unseen during training, demonstrations
from that new operator are required to provide its coarse--fine
information to the model. These demonstrations can then be used
repeatedly for multiple queries from the same operator. In contrast,
when the operator remains unchanged and only the query enters an OOD
regime, the previously available demonstrations from that operator can
be reused, and only a new coarse query needs to be generated.

We first consider parametric scalar conservation laws, where we define eight PDE families and study generalization to unseen operators, operator-parameter OOD extrapolation, and query-side OOD shifts in the initial condition. For the initial-condition OOD task, no new fine solution is required at inference time; only a new coarse query needs to be computed.
We next consider the Kuramoto--Sivashinsky equation, where we study generalization to unseen operators, operator-parameter OOD extrapolation, and cross-operator temporal OOD transfer.
We then consider the Korteweg--de Vries equation, where the model learns a temporally shifted coarse-to-fine mapping \(v_c(t)\mapsto v_f(t+\Delta t)\) with a prescribed time offset \(\Delta t\), and is evaluated on operators that are unseen during training.
Finally, we consider two-dimensional ice-sheet flow. We first compare
ICL with established multi-fidelity DeepONet baselines under a reference
multiresolution data budget, and then extend the problem to multiple
ice-sheet operators with varying rheology and basal friction to evaluate
in-context adaptation to unseen operators.


\subsection{Parametric scalar conservation laws}

\paragraph{PDE setting.}

We first consider the parametric scalar conservation
law examples studied in \cite{sun2025towards},
\[
\frac{\partial u}{\partial t}
+
q_1\frac{\partial F(u)}{\partial x}
=
q_2\frac{\partial^2u}{\partial x^2},
\qquad
x\in[0,1),\quad t\in[0,2],
\]
with periodic boundary conditions. We consider four flux functions,
$F(u)\in\{u^2/2,u^3/3,\sin(u),u\}$, each in viscous and inviscid
form, giving eight PDE families.

\paragraph{Operator setting.}

For each PDE family, an operator is defined by fixing $q_1$ and
$q_2$. We consider $q_1\in[0.5,2.0]$ for all families. For the
viscous families, $q_2\in[0.005,0.02]$, while $q_2=0$ for the
inviscid families.

For a fixed operator, different initial conditions generate different
solution samples. At the terminal time $t=2$, each demonstration is a
pair $(\mathbf v_c,\mathbf v_f)$ generated from the same initial
condition using the low- and high-fidelity numerical solvers,
respectively. For each prediction, the model is given $D$
demonstrations from the same operator together with a query
$\mathbf v_{c,q}$ generated from a different initial condition, and
predicts $\widehat{\mathbf v}_{f,q}$ without any parameter updates.

Detailed initial-condition construction, data-generation,
numerical-solver, and training settings are provided in
Appendix~\ref{app:scalar_details}.


\subsubsection{Generalization to unseen operators}
\label{sec:scalar_unseen_operators}

We first examine whether the model can generalize to operators that are
unseen during training. For each unseen testing operator, the model is
provided with a few demonstrations from that operator and predicts the
fine solution for a new query from the same operator.

\noindent\textbf{Training:}
For each of the eight PDE families, 70 training operators are generated
by sampling $q_1\in[0.5,2.0]$. For the viscous families,
$q_2\in[0.005,0.02]$ is also sampled, while $q_2=0$ for the
inviscid families. Each training example is constructed from a single
operator. The demonstrations and query correspond to different initial
conditions from that operator, and the number of demonstrations is
varied dynamically from $D=1$ to $D=5$ during training.

\noindent\textbf{Testing:}
For each PDE family, 30 unseen testing operators are generated
independently from the same $q_1$ and $q_2$ ranges used for training.
For each test case, five demonstrations are selected from one unseen
testing operator, and a query $\mathbf v_{c,q}$ is generated from a
different initial condition of that operator. The model predicts
$\widehat{\mathbf v}_{f,q}$ without any parameter updates.

\paragraph{Results.}

We evaluate the model on 1,024 fixed test cases from the unseen testing
operators across all eight PDE families. Each test case corresponds to
one query prediction from an unseen testing operator and uses five
demonstrations from that same operator. The same 1,024 test cases are
used in the subsequent scaling and context comparisons. The results are
summarized in Table~\ref{tab:sec31_unseen_results}.

\begin{table}[H]
\centering
\small
\setlength{\tabcolsep}{4.5pt}

\begin{tabular}{llcccc}
\toprule
Flux
& Regime
& \multicolumn{1}{c}{Test cases}
& \multicolumn{1}{c}{$\frac{\|\mathbf v_c-\mathbf v_f\|_2}{\|\mathbf v_f\|_2}$ (\%)}
& \multicolumn{1}{c}{$\frac{\|\widehat{\mathbf v}_f-\mathbf v_f\|_2}{\|\mathbf v_f\|_2}$ (\%)}
& \multicolumn{1}{c}{$R^2$} \\
\midrule

\textbf{Overall}
& ---
& 1024
& $\mathbf{39.33\pm26.41}$
& $\mathbf{3.66\pm2.15}$
& $\mathbf{0.9945\pm0.0044}$ \\

\midrule

$u^2/2$
& Viscous
& 111
& $28.77\pm19.61$
& $1.97\pm0.98$
& $0.9945\pm0.0047$ \\

$u^2/2$
& Inviscid
& 120
& $14.25\pm4.73$
& $4.66\pm1.63$
& $0.9896\pm0.0043$ \\

$u^3/3$
& Viscous
& 132
& $48.09\pm20.74$
& $2.16\pm0.95$
& $0.9953\pm0.0048$ \\

$u^3/3$
& Inviscid
& 114
& $13.52\pm2.11$
& $5.66\pm1.14$
& $0.9914\pm0.0031$ \\

$\sin(u)$
& Viscous
& 139
& $51.58\pm17.90$
& $2.19\pm0.87$
& $0.9968\pm0.0025$ \\

$\sin(u)$
& Inviscid
& 134
& $24.86\pm4.68$
& $5.18\pm1.16$
& $0.9942\pm0.0028$ \\

$u$
& Viscous
& 139
& $82.44\pm21.85$
& $1.64\pm0.78$
& $0.9983\pm0.0016$ \\

$u$
& Inviscid
& 135
& $40.91\pm10.89$
& $5.98\pm2.05$
& $0.9950\pm0.0035$ \\

\bottomrule
\end{tabular}

\caption{
Unseen-operator prediction across the eight scalar conservation-law
families using five same-operator coarse--fine demonstrations.
The table reports the coarse--fine relative $L^2$ discrepancy, the
ICL prediction error, and $R^2$.
}

\label{tab:sec31_unseen_results}
\end{table}

\begin{figure}[H]
    \centering
    \includegraphics[width=\textwidth]{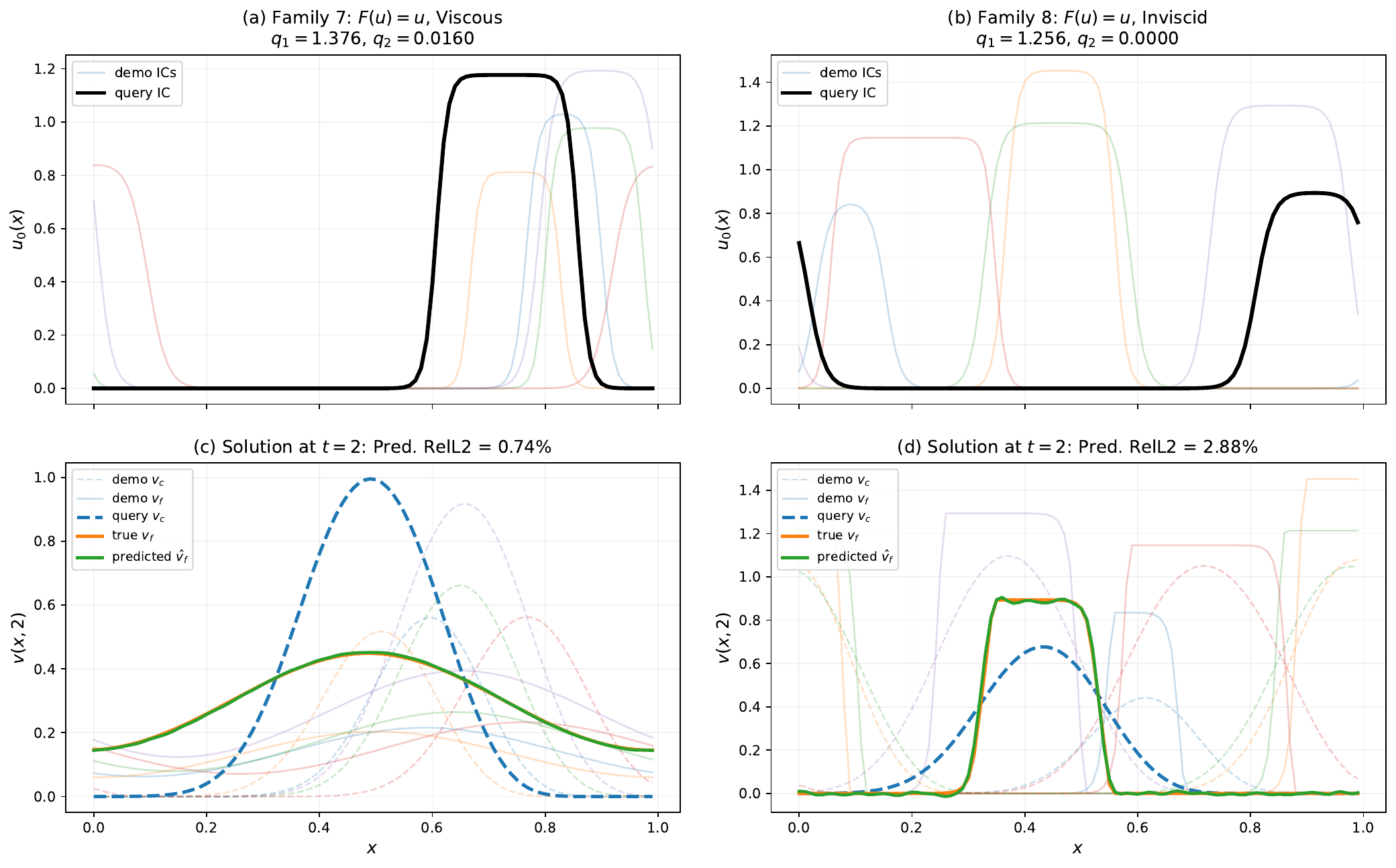}

    \caption{
Representative predictions for unseen testing scalar conservation-law
operators. The left and right columns show viscous and inviscid
linear-flux examples, respectively. The top row shows the query initial
condition in black and the demonstration initial conditions as faded
curves. The bottom row shows the solutions at $t=2$. Faded curves show
the coarse--fine demonstrations. The blue dashed curve shows the coarse
query, the orange curve shows the fine target, and the green curve shows
the ICL prediction. The relative $L^2$ prediction error is reported
above each solution panel.
    }
    \label{fig:sec31_unseen_examples}
\end{figure}

Representative predictions are shown in
Figure~\ref{fig:sec31_unseen_examples}.

\FloatBarrier

\paragraph{Results analysis.}

Across the complete 1,024-case test set, the mean coarse--fine
discrepancy is $39.33\%$, whereas the mean ICL prediction error is only
$3.66\%$. This corresponds to a reduction of approximately $90.7\%$
relative to the low-fidelity error, with a mean coefficient of
determination of $R^2=0.9945$.

Notably, the model contains only 255,364 trainable parameters while
using 560 operators for training. On average, this corresponds to only
about 456 parameters per training operator. This highlights the strong
parameter efficiency of the shared ICL model and its ability to handle
multiple operators within a single model without increasing the model
size.

The improvement is observed across all eight PDE families despite
substantial variation in the underlying fidelity gap. For example,
the viscous linear-flux family has a mean coarse--fine error of
$82.44\%$ but a prediction error of only $1.64\%$. The inviscid
families are generally more challenging, with the largest mean
prediction error, $5.98\%$, occurring for the inviscid linear-flux
family. These results show that a small set of in-context coarse--fine
demonstrations is sufficient for the trained model to adapt to
previously unseen operators.

\paragraph{Scaling and context analysis.}

\begin{figure}[H]
    \centering
    \includegraphics[width=\textwidth]{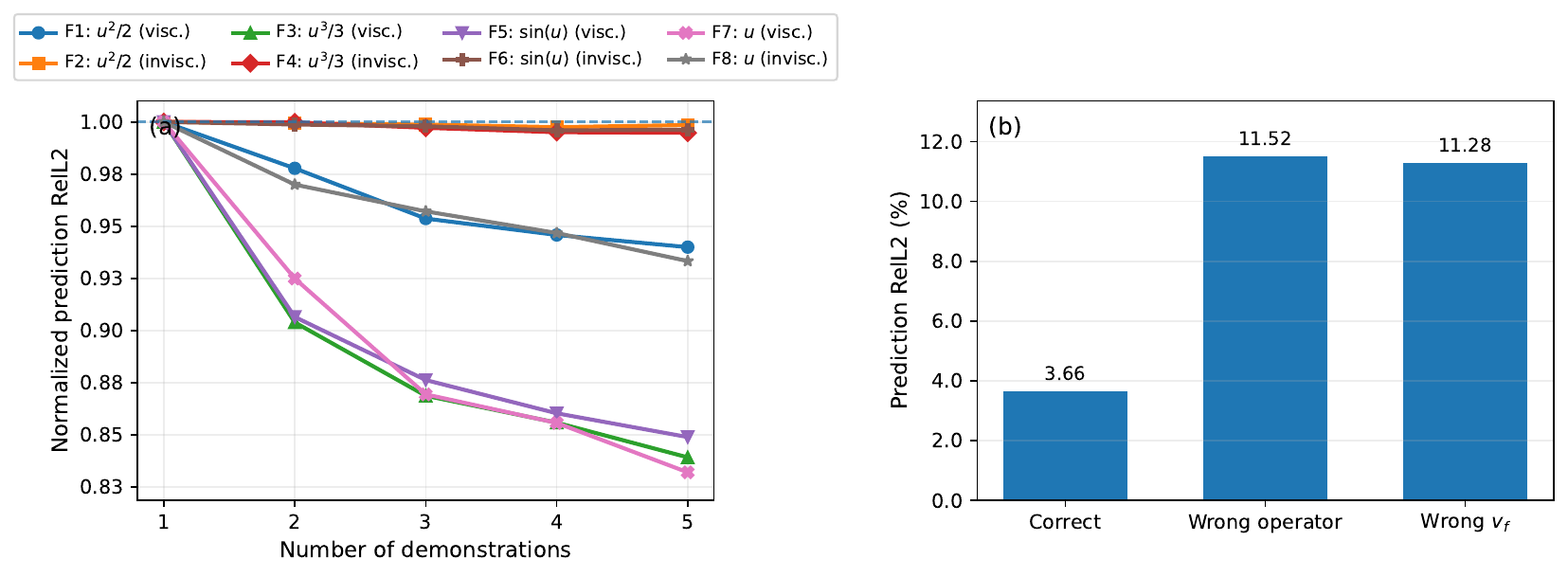}

    \caption{
Context scaling and intervention results for the scalar conservation laws.
Panel (a) shows the family-wise mean prediction error normalized by its value
with one demonstration. Panel (b) compares the mean prediction errors
under correct, wrong-operator, and wrong-$\mathbf v_f$ contexts.
The wrong-operator context replaces all demonstration pairs with
pairs from another unseen operator in the same PDE family. The wrong-$\mathbf v_f$ context retains
the demonstrated coarse solutions and replaces their fine counterparts.
All evaluations use the same trained model and fixed test queries
without parameter updates.
    }
    \label{fig:sec31_scaling_context}
\end{figure}

We next examine how prediction changes with the number and consistency
of the demonstrations. The same trained model and fixed test queries
are used in all evaluations. The results are summarized in
Figure~\ref{fig:sec31_scaling_context}.

Increasing the number of demonstrations gives a modest but consistent
improvement. The overall mean prediction error decreases from $3.88\%$
at $D=1$ to $3.77\%$, $3.71\%$, $3.68\%$, and $3.66\%$ for
$D=2$ to $D=5$, respectively. All eight PDE families have lower error
at $D=5$ than at $D=1$.

The prediction also depends strongly on the demonstration content.
With the correct five demonstrations, the mean error is $3.66\%$.
Replacing all demonstrations with those from another unseen testing
operator in the same PDE family increases the error to $11.52\%$,
while replacing only their fine solutions with those from
another unseen testing operator in the same PDE family increases it to
$11.28\%$. Both interventions increase the error by more than a factor
of three, showing that the model uses the
operator-specific coarse--fine information provided by the
demonstrations.


\subsubsection{Out-of-distribution generalization in the operator parameter}

We next examine whether the model can generalize to unseen operators
whose $q_1$ values lie outside the range used during training. In
Section~\ref{sec:scalar_unseen_operators}, the unseen testing operators are generated using the
base range $q_1\in[0.50,2.00]$. Here, the testing operators are instead
generated from the lower-OOD range $q_1\in[0.45,0.50)$ and the
upper-OOD range $q_1\in(2.00,2.20]$.

\noindent\textbf{Training:}
The training setting is the same as in Section~\ref{sec:scalar_unseen_operators}. For each of the
eight PDE families, 70 training operators are generated by sampling
$q_1\in[0.5,2.0]$. For the viscous families,
$q_2\in[0.005,0.02]$ is also sampled, while $q_2=0$ for the
inviscid families. Each training example is constructed from a single
operator. The demonstrations and query correspond to different initial
conditions from that operator, and the number of demonstrations is
varied dynamically from $D=1$ to $D=5$ during training.

\noindent\textbf{Testing:}
The OOD testing operators are generated using $q_1$ values outside the
training range. We consider a lower-OOD region,
$q_1\in[0.45,0.50)$, and an upper-OOD region,
$q_1\in(2.00,2.20]$. For the viscous families,
$q_2\in[0.005,0.02]$ is sampled as in training, while $q_2=0$ for
the inviscid families.

For each OOD test case, five demonstrations are selected from one
unseen OOD testing operator, and a query $\mathbf v_{c,q}$ is
generated from another initial condition of that operator. The model
predicts $\widehat{\mathbf v}_{f,q}$ without any parameter updates.

\paragraph{Results.}

Table~\ref{tab:sec31_operator_ood} compares unseen testing operators
from the base $q_1$ range with unseen testing operators from the lower-
and upper-OOD $q_1$ regions. Each test case corresponds to one query
prediction from an unseen testing operator in the corresponding range.

\begin{table}[H]
\centering
\small
\setlength{\tabcolsep}{4.7pt}

\begin{tabular}{lccccc}
\toprule
Setting
& $q_1$ range
& \multicolumn{1}{c}{Test cases}
& \multicolumn{1}{c}{$\frac{\|\mathbf v_c-\mathbf v_f\|_2}{\|\mathbf v_f\|_2}$ (\%)}
& \multicolumn{1}{c}{$\frac{\|\widehat{\mathbf v}_f-\mathbf v_f\|_2}{\|\mathbf v_f\|_2}$ (\%)}
& \multicolumn{1}{c}{$R^2$} \\
\midrule

Base range
& $[0.50,2.00]$
& 1024
& $39.33\pm26.41$
& $3.66\pm2.15$
& $0.9945\pm0.0044$ \\

Lower OOD
& $[0.45,0.50)$
& 498
& $45.52\pm32.41$
& $5.01\pm2.79$
& $0.9929\pm0.0054$ \\

Upper OOD
& $(2.00,2.20]$
& 526
& $31.95\pm22.92$
& $3.80\pm2.34$
& $0.9910\pm0.0090$ \\

\midrule

\textbf{Overall OOD}
& $[0.45,0.50)\cup(2.00,2.20]$
& 1024
& $\mathbf{38.55\pm28.75}$
& $\mathbf{4.39\pm2.64}$
& $\mathbf{0.9919\pm0.0075}$ \\

\bottomrule
\end{tabular}

\caption{
Operator-parameter OOD prediction for $q_1$. The base-range results
correspond to unseen operators with $q_1\in[0.50,2.00]$, while the
lower- and upper-OOD results use operators outside this range. Each
prediction uses five same-operator coarse--fine demonstrations. The
table reports the coarse--fine relative $L^2$ discrepancy, the ICL
prediction error, and $R^2$, with results shown as mean $\pm$ standard
deviation.
}
\label{tab:sec31_operator_ood}
\end{table}

\FloatBarrier

\paragraph{Results analysis.}

Across the complete OOD test set, the mean coarse--fine discrepancy is
$38.55\%$, while the mean ICL prediction error is $4.39\%$. This
corresponds to a reduction of approximately $88.6\%$ relative to the
low-fidelity error, with a mean $R^2$ of $0.9919$. Compared with the
$3.66\%$ prediction error for unseen testing operators in the base
$q_1$ range, the error increases only moderately when $q_1$ moves
outside the range used during training.

The lower-OOD region is more challenging, with a mean prediction error
of $5.01\%$, while the upper-OOD region gives $3.80\%$. The
coarse--fine discrepancy is $45.52\%$ in the lower-OOD region and
$31.95\%$ in the upper-OOD region. Overall, the model maintains
accurate prediction for unseen operators even when the operator
parameter $q_1$ lies outside the range used during training.


\subsubsection{Cross-operator initial-condition-amplitude out-of-distribution transfer}

In this section, we study whether OOD information can transfer across
operators. We consider two groups of operators, Group A and Group B.
The base-amplitude range is $A\in[0.80,1.50]$, while the near- and
farther-OOD amplitude ranges are $A\in(1.50,1.57]$ and
$A\in(1.57,1.65]$, respectively.

Importantly, this OOD task does not require any new high-fidelity
simulation at test time. The model reuses the previously available
base-amplitude coarse--fine demonstrations from the same Group B
operator and requires only a new coarse query
$\mathbf v_{c,q}$ in the OOD amplitude ranges. The corresponding
fine solution $\mathbf v_{f,q}$ is used only as ground truth
for numerical evaluation. The training and evaluation allocation is
summarized in Table~\ref{tab:sec31_amplitude_partition}.

\begin{table}[H]
\centering
\small
\setlength{\tabcolsep}{13pt}
\renewcommand{\arraystretch}{1.45}

\begin{tabular}{lcccc}
\toprule
&
\multicolumn{2}{c}{\textbf{Training}}
&
\multicolumn{2}{c}{\textbf{Evaluation}}
\\
\cmidrule(lr){2-3}
\cmidrule(lr){4-5}

&
\textbf{Group A}
&
\textbf{Group B}
&
\textbf{Group A}
&
\textbf{Group B}
\\
\midrule

\shortstack[l]{
\textbf{Base amplitude}\\
$A\in[0.80,1.50]$
}
&
$\checkmark$
&
$\checkmark$
&
--
&
--
\\

\shortstack[l]{
\textbf{Near-OOD amplitude}\\
$A\in(1.50,1.57]$
}
&
$\checkmark$
&
--
&
--
&
$\checkmark$
\\

\shortstack[l]{
\textbf{Farther-OOD amplitude}\\
$A\in(1.57,1.65]$
}
&
$\checkmark$
&
--
&
--
&
$\checkmark$
\\

\bottomrule
\end{tabular}

\caption{
Training and evaluation setup for the initial-condition amplitude OOD
experiment. Group A is trained on base-, near-OOD, and farther-OOD
amplitudes, while Group B is trained only on the base-amplitude range.
Group B is then evaluated on the near- and farther-OOD amplitude ranges
using its base-amplitude coarse--fine pairs as demonstrations.
}
\label{tab:sec31_amplitude_partition}
\end{table}

\noindent\textbf{Training:}
For each of the eight PDE families, we use 200 operators, with 100
assigned to Group A and 100 to Group B, giving 1,600 operators in
total. The operators are generated using the same $q_1$ and $q_2$
ranges as in Section~\ref{sec:scalar_unseen_operators}. We compare two training settings while
keeping the remaining training details the same.

For the \emph{base-amplitude-only} model, both Group A and Group B are
trained only on the base-amplitude range
$A\in[0.80,1.50]$.

For the \emph{mixed-amplitude} model, Group A is trained on the
base-amplitude range $A\in[0.80,1.50]$, the near-OOD range
$A\in(1.50,1.57]$, and the farther-OOD range
$A\in(1.57,1.65]$. Group B is still trained only on the
base-amplitude range $A\in[0.80,1.50]$ and receives no OOD-amplitude
data during training.

\noindent\textbf{Testing:}
Both models are evaluated on the same 800 Group B operators across
the eight PDE families. For each operator, all 100 initial-condition
samples are used to construct near- and farther-OOD queries, giving
80,000 test cases in each OOD regime. For each paired query, only the
initial-condition amplitude changes to either the near- or
farther-OOD range.

\begin{figure}[H]
    \centering
    \includegraphics[width=0.90\textwidth]{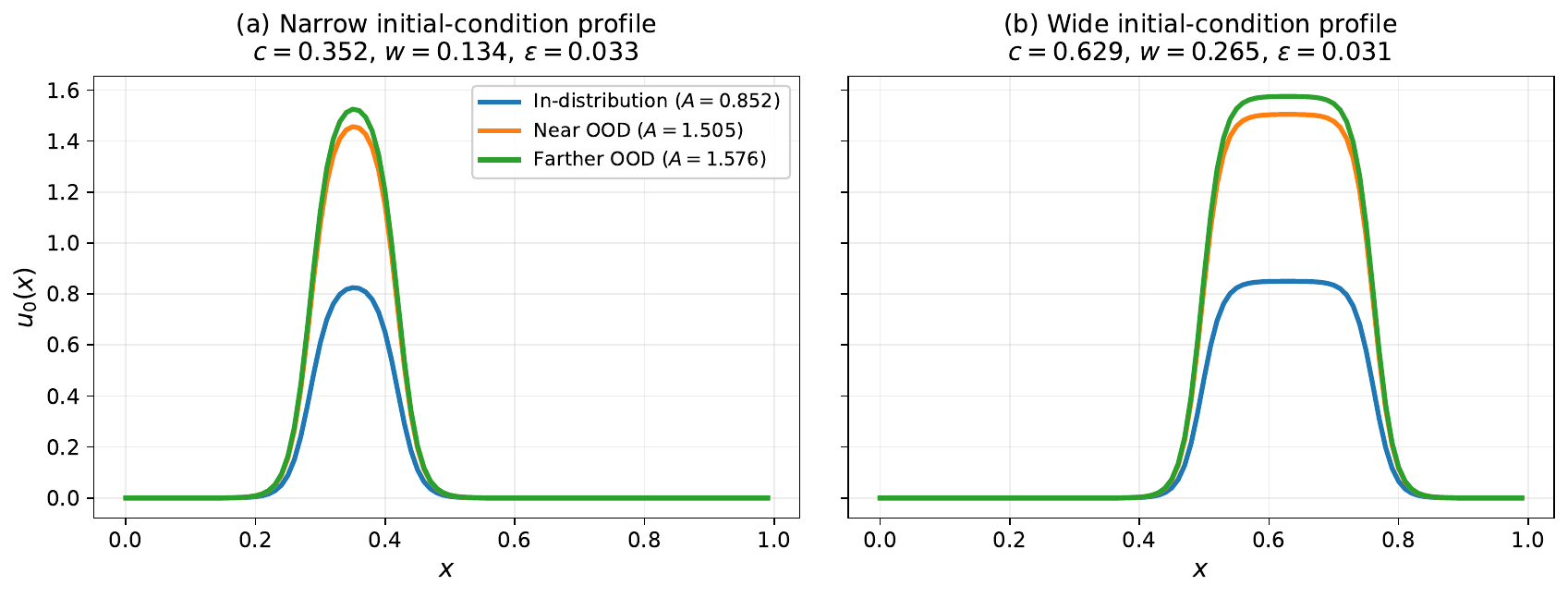}

    \caption{
Representative initial conditions for the amplitude-OOD experiment.
The left and right panels show narrow and wide profiles, respectively.
Each profile is shown at base, near-OOD, and farther-OOD amplitudes
in blue, orange, and green, respectively.
    }
    \label{fig:sec31_amplitude_examples}
\end{figure}

For each Group B query, five base-amplitude coarse--fine pairs from the
same operator are reused as demonstrations, with the query sample
excluded from the demonstration set. The base-amplitude-only and
mixed-amplitude models are evaluated using exactly the same
demonstrations for each matched query. A new coarse solution
$\mathbf v_{c,q}$ is then provided from either the near- or
farther-OOD amplitude range. Representative matched initial conditions
are shown in Figure~\ref{fig:sec31_amplitude_examples}.

\FloatBarrier

\paragraph{Results.}

Table~\ref{tab:sec31_amplitude_ood} reports the Group B test results
for the near- and farther-OOD amplitude queries. We compare two models.
The \emph{base-amplitude-only} model is trained using only
base-amplitude data from both Group A and Group B. The
\emph{mixed-amplitude} model is trained using base-, near-OOD, and
farther-OOD amplitude data from Group A and base-amplitude data from
Group B.

\begin{table}[H]
\centering
\footnotesize
\setlength{\tabcolsep}{3.0pt}
\renewcommand{\arraystretch}{1.20}

\resizebox{\textwidth}{!}{
\begin{tabular}{llcccccc}
\toprule

Query regime
& $A$ range
& Test cases
& \multicolumn{3}{c}{Relative $L^2$ error (\%)}
& \multicolumn{2}{c}{$R^2$}
\\

\cmidrule(lr){4-6}
\cmidrule(lr){7-8}

&
&
&
LF baseline
&
Base-amplitude-only model
&
Mixed-amplitude model
&
Base-amplitude-only model
&
Mixed-amplitude model
\\[-1pt]

&
&
&
$\displaystyle
\frac{\|\mathbf v_c-\mathbf v_f\|_2}{\|\mathbf v_f\|_2}$
&
$\displaystyle
\frac{\|\widehat{\mathbf v}_f-\mathbf v_f\|_2}{\|\mathbf v_f\|_2}$
&
$\displaystyle
\frac{\|\widehat{\mathbf v}_f-\mathbf v_f\|_2}{\|\mathbf v_f\|_2}$
&
&
\\

\midrule

Near OOD
& $(1.50,1.57]$
& 80,000
& $31.66\pm23.78$
& $3.88\pm2.37$
& $\mathbf{3.57\pm2.21}$
& $0.9919\pm0.0075$
& $\mathbf{0.9932\pm0.0058}$
\\

Farther OOD
& $(1.57,1.65]$
& 80,000
& $30.95\pm23.76$
& $4.18\pm2.74$
& $\mathbf{3.65\pm2.30}$
& $0.9901\pm0.0099$
& $\mathbf{0.9925\pm0.0067}$
\\

\bottomrule
\end{tabular}
}

\caption{
Group B amplitude-OOD prediction results. Each row contains 80,000
queries from 800 Group B operators. The base-amplitude-only model is
trained with base-amplitude data from both groups, while the
mixed-amplitude model additionally uses near- and farther-OOD amplitude
data from Group A. Both models use the same previously available
Group B base-amplitude demonstrations when predicting the Group B OOD
queries. 
}
\label{tab:sec31_amplitude_ood}

\end{table}

\FloatBarrier

\paragraph{Results analysis.}

For the Group B near-OOD queries, the low-fidelity baseline has a mean
error of $31.66\%$, while the base-amplitude-only model reduces the
error to $3.88\%$. The mixed-amplitude model further reduces the error
to $3.57\%$, corresponding to a relative improvement of approximately
$7.8\%$ over the base-amplitude-only model. For the farther-OOD
queries, the corresponding errors are $30.95\%$, $4.18\%$, and
$3.65\%$, giving a larger relative improvement of approximately
$12.6\%$. The mean $R^2$ also increases from $0.9919$ to $0.9932$
for near OOD and from $0.9901$ to $0.9925$ for farther OOD.

Since the Group B near- and farther-OOD data are not used during
training, the improvement shows that OOD amplitude information learned
from Group A can transfer to Group B and improve its OOD prediction.

Importantly, the model parameter count remains unchanged at 255,364
while using 1,600 operators for training. On average, this corresponds
to only about 160 parameters per training operator, showing that a large
number of operators can be handled without increasing the model size.

In practice, the amount and coverage of available data can vary across
operators. Cross-operator transfer allows the OOD information learned
from the first group to help predict the missing OOD regime for the
second group.


\subsection{Kuramoto--Sivashinsky equation}

\paragraph{PDE setting.}

We next consider the Kuramoto--Sivashinsky (KS)
equation
\[
\frac{\partial u}{\partial t}
+
u\frac{\partial u}{\partial x}
+
\frac{\partial^2u}{\partial x^2}
+
\beta\frac{\partial^4u}{\partial x^4}
=
0,
\qquad
x\in[0,2\pi),\quad t\in[0,3.5],
\]
with periodic boundary conditions.

\paragraph{Operator setting.}

We define one operator by fixing the parameter $\beta$. For the
base-range training and unseen-operator evaluation, we use
$\beta\in[0.070,0.085]$. For each operator, 20 independent solution
trajectories are generated from different initial conditions, and
multiple solution samples are obtained by evaluating each trajectory at
different times.

In this setting, each demonstration is a pair
$(\mathbf v_c(t),\mathbf v_f(t))$ evaluated at the same time and
generated using the low- and high-fidelity numerical solvers,
respectively. The multi-fidelity operator therefore maps
$\mathbf v_c(t)$ to $\mathbf v_f(t)$ at the same physical time.

For each prediction, the model is given $D$ demonstrations from the
same operator and the same trajectory at different times together with
a query $\mathbf v_{c,q}$ evaluated at another time of that trajectory.
The spatial coordinate and time associated with each demonstration and
query are also provided to the model, while the operator parameter
$\beta$ is not provided explicitly. The model predicts
$\widehat{\mathbf v}_{f,q}$ at the same query time without any parameter
updates.

For the unseen-operator experiment in Section~\ref{sec:ks_unseen_operators}, only solution states from the
late-time interval $t\in[2.50,3.50]$ are used to construct the
demonstrations and queries.
Detailed initial-condition construction, data-generation,
numerical-solver, and training settings are provided in
Appendix~\ref{app:ks_details}.


\subsubsection{Generalization to unseen operators}
\label{sec:ks_unseen_operators}

We examine whether the model can generalize to KS operators that are
unseen during training. The demonstrations and queries are constructed
from individual operator trajectories using distinct time points in
$t\in[2.50,3.50]$.

\noindent\textbf{Training:}
We generate 200 KS operators with $\beta$ distributed over
$[0.070,0.085]$. Of these, 140 operators are used for training and the
remaining 60 operators are reserved as unseen testing operators. Each
training example is constructed from a single training operator and a
single trajectory of that operator. The demonstrations and query are
taken from distinct eligible time points of the same trajectory, and
the number of demonstrations is varied dynamically from $D=1$ to
$D=5$ during training.

\noindent\textbf{Testing:}
The final evaluation uses 960 fixed test cases from the 60 unseen
testing operators, with 16 cases constructed for each operator. For
each test case, five coarse--fine demonstrations and one coarse query
are selected from a single trajectory of one unseen testing operator
using distinct time points in $t\in[2.50,3.50]$. The model predicts
$\widehat{\mathbf v}_{f,q}$ without any parameter updates.

\paragraph{Results.}

We evaluate the model on the 960 fixed test cases from the 60 unseen
testing KS operators. Each test case corresponds to one query
prediction and uses five demonstrations from the same operator and the
same trajectory. The results are summarized in
Table~\ref{tab:sec35_ks_unseen_results}.

\begin{table}[H]
\centering
\small
\setlength{\tabcolsep}{4.5pt}

\begin{tabular}{lcccc}
\toprule
Setting
& \multicolumn{1}{c}{Test cases}
& \multicolumn{1}{c}{$\frac{\|\mathbf v_c-\mathbf v_f\|_2}{\|\mathbf v_f\|_2}$ (\%)}
& \multicolumn{1}{c}{$\frac{\|\widehat{\mathbf v}_f-\mathbf v_f\|_2}{\|\mathbf v_f\|_2}$ (\%)}
& \multicolumn{1}{c}{$R^2$} \\
\midrule

Unseen operators
& 960
& $29.23 \pm 35.25$
& $4.39 \pm 5.96$
& $0.994$ \\

\bottomrule
\end{tabular}

\caption{
Prediction accuracy on 60 unseen testing KS operators using five
same-operator, same-trajectory coarse--fine demonstrations. The
coarse query, fine target, and ICL prediction correspond
to the same query time. $R^2$ is reported as the mean over the 960 fixed test cases.
}
\label{tab:sec35_ks_unseen_results}
\end{table}

\FloatBarrier

\paragraph{Results analysis.}

Across the 960 fixed test cases, the mean relative $L^2$ discrepancy
between the coarse query and the fine target is
$29.23\%$, whereas the mean ICL prediction error is reduced to
$4.39\%$. This corresponds to an $84.98\%$ reduction in the mean
relative error. The corresponding mean coefficient of determination is
$R^2=0.994$.

The median prediction error is $2.62\%$, and the 90th-percentile
prediction error is $8.86\%$. These results show that the shared
ICL model generalizes effectively to KS operators not
observed during training.

Because the operator parameter $\beta$ is not provided explicitly to
the model, prediction on an unseen operator relies on the coarse--fine
information supplied through the in-context demonstrations together
with the explicitly provided spatial and temporal information.

\begin{figure}[H]
\centering
\includegraphics[width=\textwidth]{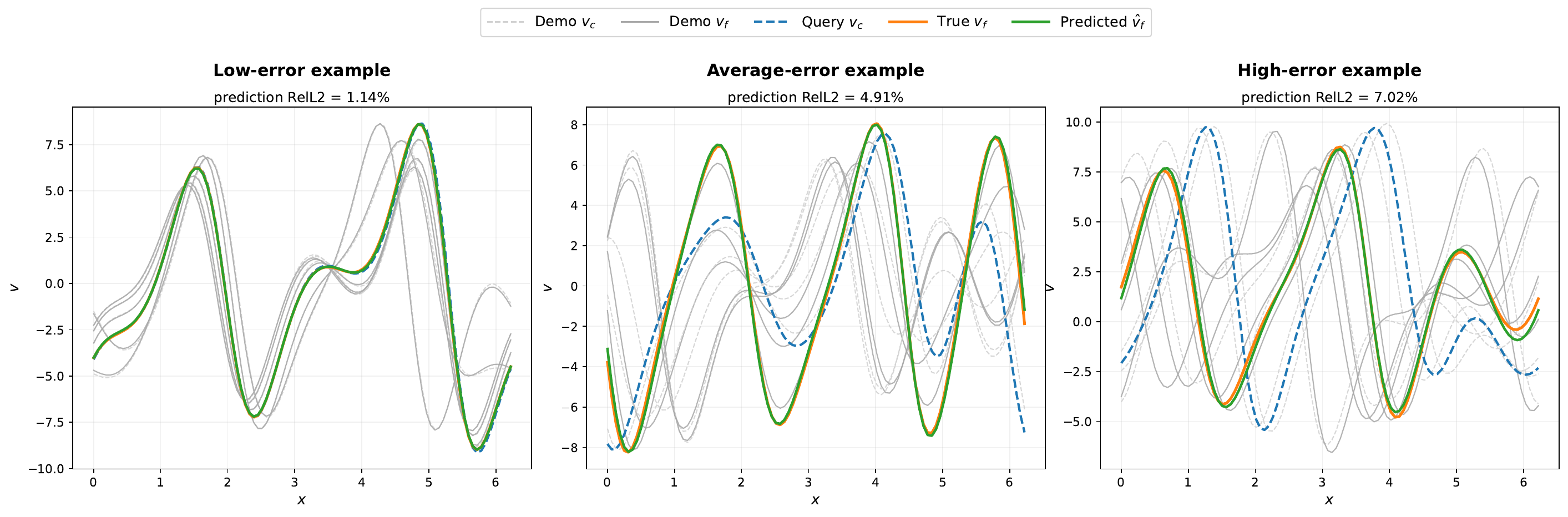}
\caption{
Representative predictions for unseen testing KS operators.
The panels show low-, average-, and high-error examples.
Faded curves show the coarse--fine demonstrations. The blue dashed
curve shows the coarse query, the orange curve shows the fine target,
and the green curve shows the ICL prediction.
The relative $L^2$ prediction error is reported above each panel.
}
\label{fig:sec35_ks_examples}
\end{figure}

Figure~\ref{fig:sec35_ks_examples} shows representative low-, average-,
and high-error predictions selected from the unseen-operator test set.
Although the coarse and fine states are evaluated at the same
physical time, the coarse query may exhibit noticeable phase,
amplitude, and local-shape differences relative to the high-fidelity
target. The model uses the supplied coarse--fine demonstrations to
correct these discrepancies and recover the high-fidelity solution at the query time.

\FloatBarrier

\paragraph{Scaling and context analysis.}

We next examine how the model responds to both the amount and the
consistency of the supplied context. The same trained model and the same
960 fixed test queries are used throughout, with no parameter updates.
For the scaling study, the demonstrations are nested so that the
context with $D$ demonstrations consists of the first $D$
demonstrations of the corresponding five-demonstration test case.

The mean prediction error decreases monotonically from $9.37\%$ at
$D=1$ to $6.63\%$, $5.51\%$, $4.85\%$, and $4.39\%$ for
$D=2,\ldots,5$, respectively. After normalization by the corresponding
$D=1$ error, the mean errors are $0.707$, $0.588$, $0.517$, and $0.469$
for $D=2,\ldots,5$, respectively. Thus, increasing the number of valid
demonstrations produces a systematic improvement in prediction
accuracy.

The prediction also depends strongly on the consistency of the
demonstration content. With the correct five demonstrations, the mean
prediction error is $4.39\%$. Replacing the demonstrations by coherent
coarse--fine pairs from another unseen operator increases the mean error
to $52.10\%$. Retaining the original demonstrated coarse states while
replacing only the demonstrated fine states increases the mean
error further to $127.28\%$.

\begin{figure}[H]
\centering
\includegraphics[width=\textwidth]{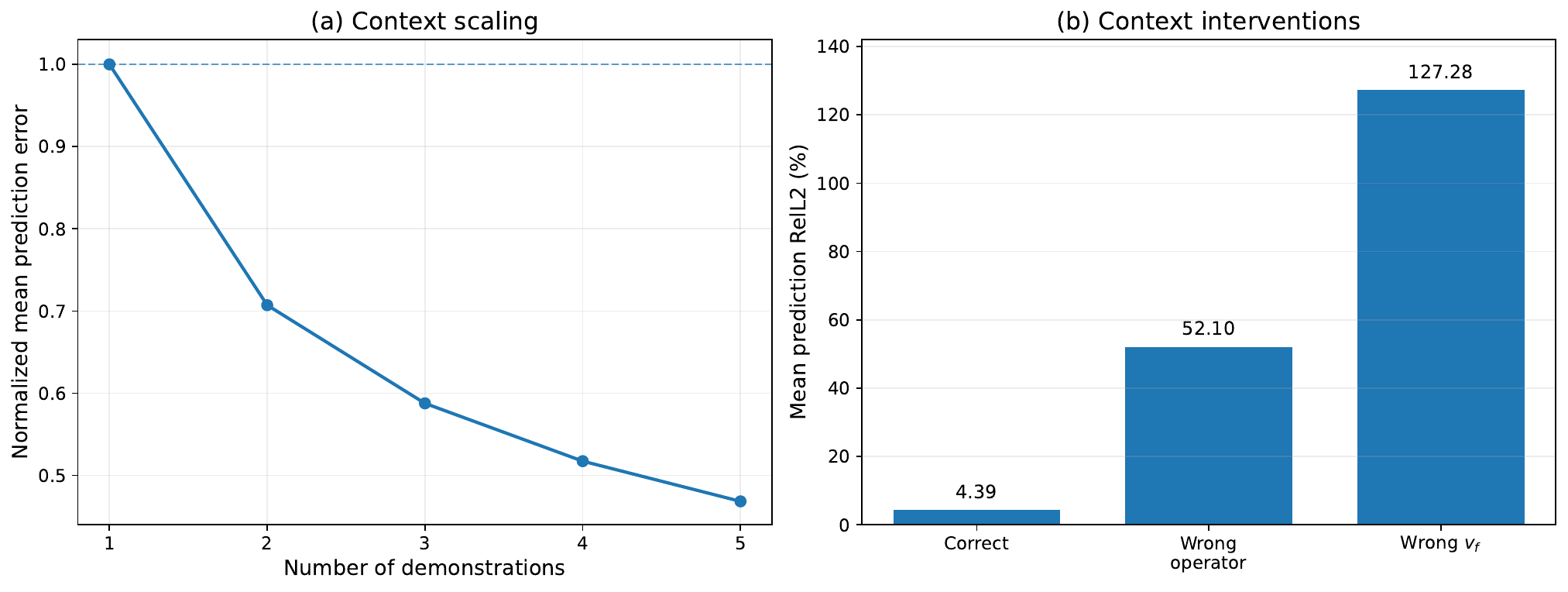}
\caption{
Context scaling and intervention results for KS.
Panel (a) shows the mean prediction error normalized by its value
with one demonstration. Panel (b) compares the mean prediction errors
under correct, wrong-operator, and wrong-$\mathbf v_f$ contexts.
The wrong-operator context replaces all demonstration pairs with
pairs from another unseen operator. The wrong-$\mathbf v_f$ context retains
the demonstrated coarse solutions and replaces their fine counterparts.
All evaluations use the same trained model and fixed test queries
without parameter updates.
}
\label{fig:sec35_ks_scaling_intervention}
\end{figure}

Figure~\ref{fig:sec35_ks_scaling_intervention} summarizes the scaling
and context-intervention results.

\FloatBarrier

These intervention results show that the model does not simply rely on
the query state alone. Instead, the prediction depends strongly on the
coarse--fine relations encoded by the demonstrations and on their
dynamic consistency with the current query trajectory. In particular,
the wrong-$\mathbf v_f$ intervention directly disrupts the demonstrated
coarse-to-fine correspondence and produces the largest degradation in
prediction accuracy. Together with the monotonic scaling behavior,
these results provide evidence that the model actively uses the
in-context demonstrations when constructing the high-fidelity
prediction.


\subsubsection{Out-of-distribution generalization in the operator parameter}
\label{sec:ks_beta_ood}

We next examine whether the model can generalize to unseen KS operators
whose $\beta$ values lie outside the range used during training. The
base range is $\beta\in[0.070,0.085]$, and the OOD operators use
$\beta\in(0.085,0.100]$.

\noindent\textbf{Training:}
The training setting is the same as in
Section~\ref{sec:ks_unseen_operators}. The model is trained using 140 operators with
$\beta\in[0.070,0.085]$, with 20 trajectories per operator. Each
training example is constructed from one operator and one trajectory,
and the number of demonstrations is varied dynamically from $D=1$ to
$D=5$. No OOD operator is used for training, normalization, or
checkpoint selection.

\noindent\textbf{Testing:}
We generate 50 unseen OOD operators with
$\beta\in(0.085,0.100]$. The evaluation contains 512 fixed test cases.
For each prediction, the five coarse--fine demonstrations and the query
come from the same OOD operator and the same trajectory. The trained
model predicts the corresponding fine solution without any parameter
updates.

\paragraph{Results.}

Table~\ref{tab:sec35_ks_operator_ood} compares the base-range unseen
operators with the operator-parameter OOD operators.

\begin{table}[H]
\centering
\small
\setlength{\tabcolsep}{5pt}

\begin{tabular}{lccccc}
\toprule
Setting
& $\beta$ range
& \multicolumn{1}{c}{Test cases}
& \multicolumn{1}{c}{$\frac{\|\mathbf v_c-\mathbf v_f\|_2}{\|\mathbf v_f\|_2}$ (\%)}
& \multicolumn{1}{c}{$\frac{\|\widehat{\mathbf v}_f-\mathbf v_f\|_2}{\|\mathbf v_f\|_2}$ (\%)}
& \multicolumn{1}{c}{$R^2$} \\
\midrule

Base range
& $[0.070,0.085]$
& 960
& $29.23\pm35.25$
& $4.39\pm5.96$
& $0.9942\pm0.0376$ \\

OOD
& $(0.085,0.100]$
& 512
& $17.70\pm25.52$
& $\mathbf{4.63\pm4.72}$
& $\mathbf{0.9953\pm0.0208}$ \\

\bottomrule
\end{tabular}

\caption{
Operator-parameter OOD prediction for the KS equation. The base-range
results correspond to unseen operators with
$\beta\in[0.070,0.085]$, while the OOD results use unseen operators
with $\beta\in(0.085,0.100]$. Each OOD prediction uses five
same-operator, same-trajectory coarse--fine demonstrations.
}
\label{tab:sec35_ks_operator_ood}
\end{table}

\FloatBarrier

\paragraph{Results analysis.}

For $\beta\in(0.085,0.100]$, the mean ICL prediction error is
$4.63\%$, compared with $4.39\%$ for unseen operators in the base
range. The mean coefficient of determination remains high at
$R^2=0.9953$. The model therefore maintains prediction accuracy close
to the base-range level when the operator parameter moves beyond the
range used during training.


\subsubsection{Cross-operator temporal out-of-distribution transfer}

We next examine whether temporal OOD information learned from one group
of KS operators can improve prediction for another group that has not
observed the same OOD temporal regime. We consider two temporal OOD
intervals, $[3.002,3.100]$ and $[3.005,3.250]$. For each interval,
we compare two matched models: a \emph{base-time-only} model and a
\emph{mixed-time} model.

For each experiment, the 200 KS operators are divided into two fixed
groups of 100 operators, denoted Group A and Group B. Each operator
contains 20 independent trajectories. Both experiments use 100
base-time points from $t=2.505$ to $t=3.000$. Each temporal OOD set
contains 50 time points: $t=3.002,3.004,\ldots,3.100$ in the first
experiment and $t=3.005,3.010,\ldots,3.250$ in the second.

For the \emph{base-time-only} model, both Group A and Group B are
trained using only the base-time points. For the \emph{mixed-time}
model, Group A is trained using both its 100 base-time points and 50 OOD
time points, while Group B is trained using only its 100 base-time
points. Thus, temporal OOD information is introduced only through
Group A.

As in the standard KS experiment, each training example is constructed
from one operator and one trajectory, and the demonstrations and query
are selected from distinct eligible time points of that trajectory.
During training, the number of demonstrations is varied dynamically
from $D=1$ to $D=5$. For the mixed-time model, the Group A eligible
time pool contains 100 base-time and 50 OOD points, so one third of its
eligible time points lie in the OOD interval. Group B OOD states are
not used during training.

\noindent\textbf{Testing:}
Both models are evaluated on the same Group B temporal OOD queries. For
each of the 100 Group B operators, all 20 trajectories and all 50 OOD
query times are evaluated, giving
$100\times20\times50=100{,}000$ paired test cases for each temporal
OOD interval.

Each test case uses five coarse--fine demonstrations from the base-time
interval of the same Group B operator and the same trajectory as the
query. The two models are evaluated using the same base-time
demonstration context for each paired query. A coarse query
$\mathbf v_{c,q}$ from the OOD interval is then provided to the model,
while the corresponding fine solution is used only for
evaluation.

The Group A/Group B allocation for the mixed-time setting is summarized
in Table~\ref{tab:sec35_ks_temporal_partition}. The base-time-only model
uses the same operator partition, with both groups trained only on the
base-time set.

\begin{table}[H]
\centering
\small
\setlength{\tabcolsep}{13pt}
\renewcommand{\arraystretch}{1.45}

\begin{tabular}{lcccc}
\toprule
&
\multicolumn{2}{c}{\textbf{Training}}
&
\multicolumn{2}{c}{\textbf{Evaluation}}
\\
\cmidrule(lr){2-3}
\cmidrule(lr){4-5}

&
\textbf{Group A}
&
\textbf{Group B}
&
\textbf{Group A}
&
\textbf{Group B}
\\
\midrule

\shortstack[l]{
\textbf{Base time}\\
$t\in[2.505,3.000]$
}
&
$\checkmark$
&
$\checkmark$
&
--
&
--
\\

\shortstack[l]{
\textbf{Temporal OOD}\\
$t\in[3.002,3.100]$ or\\
$t\in[3.005,3.250]$
}
&
$\checkmark$
&
--
&
--
&
$\checkmark$
\\

\bottomrule
\end{tabular}

\caption{
Data allocation for each mixed-time KS temporal OOD experiment. Group A
is trained on both base-time and temporal OOD data, while Group B is
trained only on base-time data. Group B is then evaluated at temporal
OOD times using its base-time coarse--fine demonstrations as context.
}
\label{tab:sec35_ks_temporal_partition}
\end{table}

\paragraph{Results.}

Table~\ref{tab:sec35_ks_temporal_ood} reports the Group B temporal OOD
results. The \emph{base-time-only} model is trained using only base-time
data from both groups, while the \emph{mixed-time} model additionally
uses temporal OOD data from Group A.

\begin{table}[H]
\centering
\footnotesize
\setlength{\tabcolsep}{3.0pt}
\renewcommand{\arraystretch}{1.20}

\resizebox{\textwidth}{!}{
\begin{tabular}{llcccccc}
\toprule

Query regime
& $t$ range
& Test cases
& \multicolumn{3}{c}{Relative $L^2$ error (\%)}
& \multicolumn{2}{c}{$R^2$}
\\

\cmidrule(lr){4-6}
\cmidrule(lr){7-8}

&
&
&
LF baseline
&
Base-time-only model
&
Mixed-time model
&
Base-time-only model
&
Mixed-time model
\\[-1pt]

&
&
&
$\displaystyle
\frac{\|\mathbf v_c-\mathbf v_f\|_2}{\|\mathbf v_f\|_2}$
&
$\displaystyle
\frac{\|\widehat{\mathbf v}_f-\mathbf v_f\|_2}{\|\mathbf v_f\|_2}$
&
$\displaystyle
\frac{\|\widehat{\mathbf v}_f-\mathbf v_f\|_2}{\|\mathbf v_f\|_2}$
&
&
\\

\midrule

Temporal OOD
& $[3.002,3.100]$
& 100,000
& $28.77\pm34.02$
& $2.90\pm2.90$
& $\mathbf{2.70\pm2.67}$
& $0.9982$
& $\mathbf{0.9984}$
\\

Temporal OOD
& $[3.005,3.250]$
& 100,000
& $30.52\pm35.34$
& $6.62\pm9.02$
& $\mathbf{5.08\pm7.02}$
& $0.9865$
& $\mathbf{0.9919}$
\\

\bottomrule
\end{tabular}
}

\caption{
Cross-operator temporal OOD prediction on 100 Group B KS operators for
two evaluation intervals. Each row contains 100,000 paired test cases
from 20 trajectories and 50 OOD query times per operator.
The base-time-only
model is trained on base-time data from Groups A and B, while the
mixed-time model additionally uses temporal OOD data from Group A.
Group B OOD states are not used during training. Each prediction uses
five same-operator, same-trajectory base-time coarse--fine
demonstrations. Relative $L^2$ errors are reported as mean $\pm$
standard deviation, and $R^2$ is reported as the mean.
}
\label{tab:sec35_ks_temporal_ood}
\end{table}

\FloatBarrier

\paragraph{Results analysis.}

For query times in $[3.002,3.100]$, the low-fidelity baseline has a
mean relative $L^2$ error of $28.77\%$. The mixed-time model reduces
the mean prediction error from $2.90\%$ for the base-time-only model
to $2.70\%$, a relative improvement of $6.72\%$. The paired win rate
is $57.71\%$, and the mean $R^2$ increases from $0.9982$ to $0.9984$.

For query times in $[3.005,3.250]$, the low-fidelity baseline has a
mean relative $L^2$ error of $30.52\%$. The mixed-time model reduces
the mean prediction error from $6.62\%$ to $5.08\%$, a relative
improvement of $23.33\%$. The paired win rate is $69.62\%$, and the
mean $R^2$ increases from $0.9865$ to $0.9919$. Both experiments show
a transfer benefit, with a larger improvement in the experiment
covering later query times.

These results show cross-operator temporal transfer in the KS setting.
Temporal OOD information introduced only through Group A improves
prediction for Group B, even though Group B is trained only on
base-time states. The previously available Group B base-time
coarse--fine demonstrations can therefore be reused for OOD prediction,
so a new temporal OOD query requires only its coarse solution as
model input.


\subsection{Korteweg--de Vries equation}

\paragraph{PDE setting.}

We next consider the Korteweg--de Vries (KdV)
equation
\[
\frac{\partial u}{\partial t}
+
u\frac{\partial u}{\partial x}
+
\epsilon\frac{\partial^3u}{\partial x^3}
=
0,
\qquad
x\in[0,2\pi],\quad t\in[0,10],
\]
with periodic boundary conditions.

\paragraph{Operator setting.}

We define one operator by fixing the parameter triple
$(\epsilon,A,\phi)$, where $\epsilon\in[0.01,0.02]$,
$A\in[0.35,0.65]$, and $\phi\in[0,2\pi)$. For each operator,
the governing equation generates an entire solution trajectory, and
multiple solution samples are obtained by evaluating this trajectory
at different times.

In this setting, each demonstration is a pair
$(\mathbf v_c(t),\mathbf v_f(t+\Delta t))$ with
$\Delta t=0.5$. The multi-fidelity operator in this experiment maps a coarse
state at time $t$ to a fine state at the later time
$t+\Delta t$.
In the notation of Section~\ref{sec:mf_operator_formulation}, the input $u=t$ is the source time:
$G_i^c(t)$ returns the coarse state at $t$, while $G_i^f(t)$
returns the corresponding fine target at $t+\Delta t$ for this
temporally shifted task.

For each prediction, the model is given $D$ demonstrations from the
same operator at different source times together with a query
$\mathbf v_{c,q}$ evaluated at another source time of that operator.
The model predicts the corresponding fine solution
$\widehat{\mathbf v}_{f,q}$ at $t_q+0.5$ without any parameter
updates.

Detailed initial-condition construction, data-generation,
numerical-solver, and training settings are provided in
Appendix~\ref{app:kdv_details}.


\subsubsection{Generalization to unseen operators}

We examine whether the model can generalize to KdV operators that are
unseen during training. The demonstrations and queries are constructed
from the operator trajectories using source times in $t\in[0,9]$,
with the corresponding fine targets evaluated at $t+0.5$.

\noindent\textbf{Training:}
We generate 140 training operators by independently sampling
$\epsilon\in[0.01,0.02]$, $A\in[0.35,0.65]$, and
$\phi\in[0,2\pi)$. Each training example is constructed from a
single operator. The demonstrations and query are taken from different
source times of that operator, and the number of demonstrations is
varied dynamically from $D=1$ to $D=5$ during training.

\noindent\textbf{Testing:}
We use 60 additional operators as unseen testing operators, generated
independently from the same parameter ranges. These operators are not
used during training. For each test case, five demonstrations are
selected from one unseen testing operator, and a coarse query
$\mathbf v_{c,q}$ is taken from another source time of that operator.
The model predicts the corresponding fine solution at
$t_q+0.5$ without any parameter updates.

\paragraph{Results.}

We evaluate the model on 960 fixed test cases from the 60 unseen
testing KdV operators, with 16 query cases per operator. Each
prediction uses five demonstrations from the same testing operator.
The same 960 fixed test cases are used in the subsequent scaling and
context comparisons.

For the following results, $\mathbf v_c$ denotes the coarse query
at source time $t_q$, while $\mathbf v_f$ denotes the corresponding
fine target at $t_q+0.5$. The results are summarized in
Table~\ref{tab:sec33_kdv_unseen_results}.

\begin{table}[H]
\centering
\small
\setlength{\tabcolsep}{4.5pt}

\begin{tabular}{lcccc}
\toprule
Setting
& \multicolumn{1}{c}{Test cases}
& \multicolumn{1}{c}{$\frac{\|\mathbf v_c-\mathbf v_f\|_2}{\|\mathbf v_f\|_2}$ (\%)}
& \multicolumn{1}{c}{$\frac{\|\widehat{\mathbf v}_f-\mathbf v_f\|_2}{\|\mathbf v_f\|_2}$ (\%)}
& \multicolumn{1}{c}{$R^2$} \\
\midrule

Unseen operators
& 960
& $20.18\pm6.77$
& $\mathbf{3.33\pm3.52}$
& $0.9977\pm0.0070$ \\

\bottomrule
\end{tabular}

\caption{
Prediction accuracy on 60 unseen testing KdV operators using five
same-operator coarse--fine demonstrations. Here $\mathbf v_c$ denotes
the coarse query at source time $t_q$, while $\mathbf v_f$ and
$\widehat{\mathbf v}_f$ denote the fine target and ICL
prediction at $t_q+0.5$, respectively. 
}
\label{tab:sec33_kdv_unseen_results}
\end{table}

\begin{figure}[H]
    \centering
    \includegraphics[width=0.95\textwidth]
    {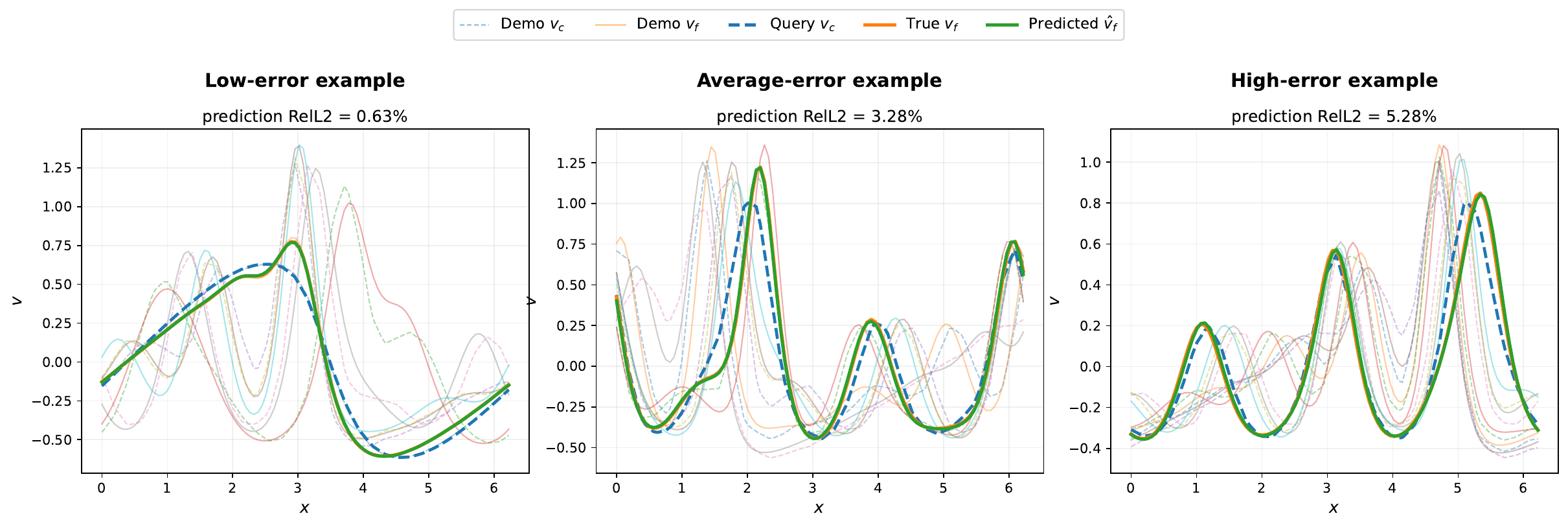}

    \caption{
Representative predictions for unseen testing KdV operators.
The panels show low-, average-, and high-error examples.
Faded curves show the coarse--fine demonstrations. The blue dashed
curve shows the coarse query, the orange curve shows the fine target,
and the green curve shows the ICL prediction.
The relative $L^2$ prediction error is reported above each panel.
The coarse query is evaluated at $t_q$, while the fine target and
prediction correspond to $t_q+0.5$.
    }
    \label{fig:sec33_kdv_unseen_example}
\end{figure}

Representative prediction results are shown in
Figure~\ref{fig:sec33_kdv_unseen_example}.

\FloatBarrier

\paragraph{Results analysis.}

Across the 960 fixed test cases, the mean relative $L^2$ discrepancy
between the coarse query at $t_q$ and the fine target at $t_q+0.5$
is $20.18\%$, whereas the mean ICL prediction error is $3.33\%$.
This discrepancy reflects both the difference in fidelity and the
change in the solution over time. The ICL prediction error is
approximately $83.5\%$ lower than this discrepancy, with a mean
coefficient of determination of $R^2=0.9977$.

The median prediction error is $2.34\%$, and the 90th-percentile
prediction error is $6.97\%$. These results show that the shared
ICL model generalizes effectively to KdV operators not
observed during training.

Figure~\ref{fig:sec33_kdv_unseen_example} shows three representative
predictions from unseen testing operators, each including the five
same-operator demonstrations, coarse query, fine target,
and ICL prediction.

\paragraph{Scaling and context analysis.}

\begin{figure}[H]
    \centering
    \includegraphics[width=0.90\textwidth]
    {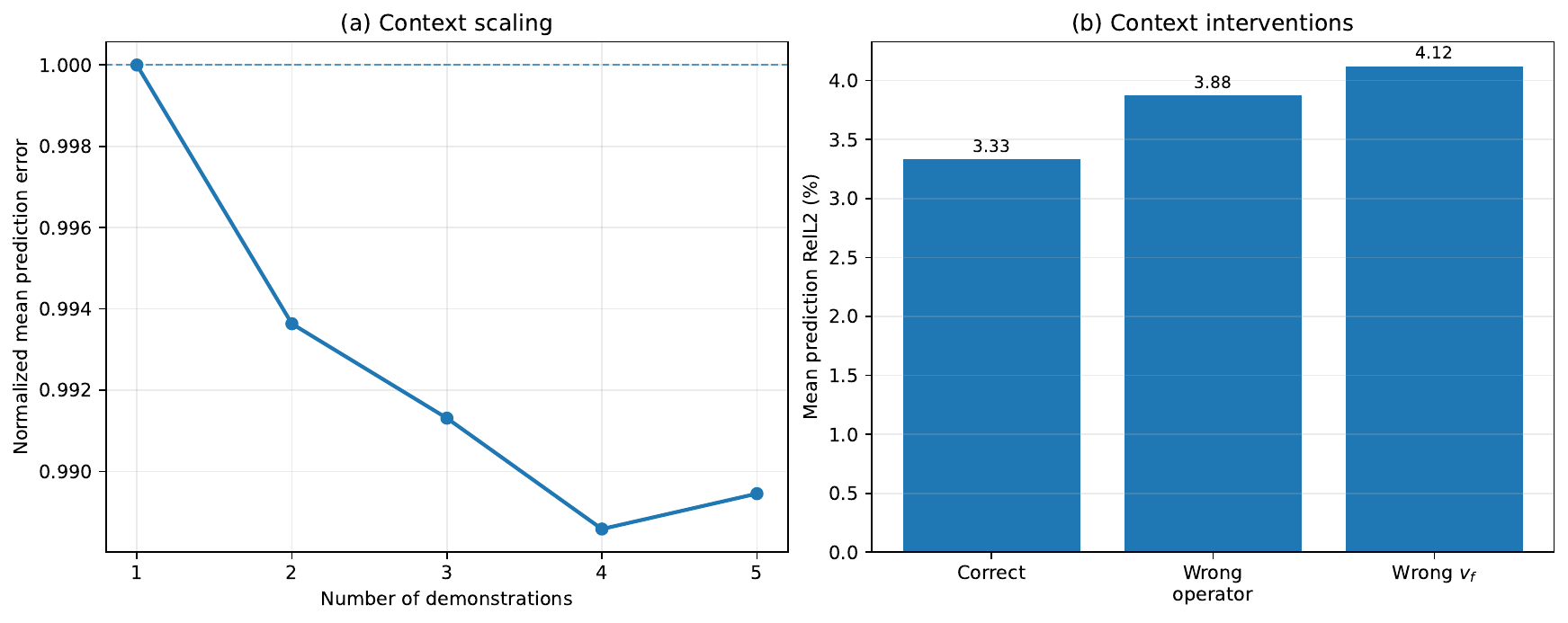}

    \caption{
Context scaling and intervention results for KdV.
Panel (a) shows the mean prediction error normalized by its value
with one demonstration. Panel (b) compares the mean prediction errors
under correct, wrong-operator, and wrong-$\mathbf v_f$ contexts.
The wrong-operator context replaces all demonstration pairs with
pairs from another unseen operator. The wrong-$\mathbf v_f$ context retains
the demonstrated coarse solutions and replaces their fine counterparts.
All evaluations use the same trained model and fixed test queries
without parameter updates.
    }
    \label{fig:sec33_kdv_scaling_context}
\end{figure}

We next examine how the model responds to both the amount and the
consistency of the supplied context. The same trained model and the
same 960 fixed test query cases are used throughout, with no parameter
updates. For the scaling study, the demonstrations are nested so that
the context with $D$ demonstrations consists of the first $D$
demonstrations of the corresponding five-demonstration test case.
The results are shown together in
Figure~\ref{fig:sec33_kdv_scaling_context}.

Increasing the number of valid demonstrations produces a modest
improvement. The overall mean prediction error decreases from
$3.37\%$ at $D=1$ to $3.34\%$, $3.34\%$, $3.33\%$, and
$3.33\%$ for $D=2,\ldots,5$, respectively.

After normalization by the corresponding $D=1$ error, the mean errors
are $0.9936$, $0.9913$, $0.9886$, and $0.9895$ for
$D=2,\ldots,5$, respectively. The improvement is not strictly
monotonic: the normalized error increases slightly from $D=4$ to $D=5$.

The context interventions produce a larger effect. With the correct
five demonstrations, the mean prediction error is $3.33\%$.
Replacing all demonstration pairs with coarse--fine pairs from another
unseen testing operator increases the error to $3.88\%$. Retaining
the original demonstrated coarse states while replacing only their
fine counterparts increases the error further to $4.12\%$.

The corresponding mean $R^2$ values are $0.9977$, $0.9971$, and
$0.9968$ for the correct, wrong-operator, and wrong-$\mathbf v_f$
contexts, respectively. Taken together, the scaling and intervention
results show that the model uses the operator-specific coarse--fine
information supplied through the in-context demonstrations.


\subsection{Two-dimensional ice-sheet flow}

\paragraph{PDE setting.}

We next consider an ice-sheet flow problem motivated by a
previously established multiresolution multi-fidelity benchmark
\cite{howard2023multifidelity,he2023hybriddeepneuraloperatorfinite}. The depth-averaged horizontal velocity $\vec{u}=(u,v)$ satisfies

\[
-\nabla\cdot(\kappa\nabla u)+\beta u
=
-\rho g H H_x,
\qquad
-\nabla\cdot(\kappa\nabla v)+\beta v
=
-\rho g H H_y,
\]

while the ice thickness evolves according to
$H_t+\nabla\cdot(H\vec{u})=0$. Here, $H$ denotes the ice thickness,
$\beta$ is the basal-friction field, and $\kappa$ represents the effective
rheological coefficient. The constitutive relations and physical
parameters follow the reference benchmark.

\paragraph{Operator setting.}

We consider two settings for this problem. The first follows the
reference multi-fidelity benchmark and uses a single fixed operator.
Different basal-friction fields are treated as different input samples
under the same reference setting, and each sample contributes one
solution at the final time $t=50$.

The second extends this benchmark to a multi-operator ICL setting. In
this extension, one operator is defined by fixing the rheological
parameter $A$ and the basal-friction field $\beta(x,y)$. Each operator
provides solution states at different times, allowing the model to use
coarse--fine demonstrations from that operator as context for new
predictions.

In the multi-operator setting, each demonstration is a pair
$(\mathbf v_c,\mathbf v_f)$ from the same operator and evaluated at
the same time. A coarse solution $\mathbf v_{c,q}$ from another time
of that operator is used as the query, and the model predicts
$\widehat{\mathbf v}_{f,q}$ without any parameter updates.

The coarse and fine solutions are generated from the same
governing equations at different spatial resolutions. Detailed
initial-condition construction, data-generation, numerical-solver,
and training settings are provided in
Appendix~\ref{app:ice_details}.


\subsubsection{Reference multi-fidelity benchmark}

We first consider a single-operator multi-fidelity setting motivated
by the reference multi-fidelity benchmark
\cite{howard2023multifidelity}. The purpose of this experiment is to
compare how effectively different methods use the reference budget of
20 high-fidelity velocity samples. Different basal-friction samples
are treated as different input samples under the same reference
setting, and each sample contributes one solution at the final time
$t=50$.

\noindent\textbf{Training:}
We use 100 low-fidelity samples and 20 high-fidelity velocity samples,
following the data budget of the reference multi-fidelity benchmark.
The 20 high-fidelity velocity samples are paired with their
corresponding coarse solutions.

For ICL, only these 20 paired coarse--fine samples are used for
training. Each training example is constructed from different paired
samples under the same reference setting, and the number of
demonstrations is sampled independently for each example from
$D=1$ to $D=5$ during training. The application-specific
information $a(u)$ contains the spatial coordinates $XY$, the
basal-friction field $\beta$, and the ice thickness $H$.
These quantities are represented by $n_a=3$ auxiliary rows,
so that $\boldsymbol{\alpha}_s$ expands into
$\mathbf{XY}_s,\boldsymbol{\beta}_s,\mathbf H_s$ for each
demonstration or query sample $s$. Each demonstration contains
$[\boldsymbol{\alpha}_{d_j},\mathbf v_{c,d_j},\mathbf v_{f,d_j}]$,
while the query contains
$[\boldsymbol{\alpha}_q,\mathbf v_{c,q}]$.
Because the velocity solution depends on the ice thickness,
including $H$ provides additional physical information for the
coarse-to-fine prediction.

We also consider the composite and non-composite multi-fidelity
DeepONet formulations introduced in the reference benchmark
\cite{howard2023multifidelity}. The composite model
uses all 100 low-fidelity samples together with the 20 high-fidelity
velocity samples. The non-composite model directly uses the available
coarse solutions in the coarse-to-fine model and uses the 20
low-fidelity samples paired with the 20 high-fidelity velocity
training samples. Thus, all three learned models use the same
20-sample high-fidelity velocity training budget, while the composite
model uses an additional 80 low-fidelity training samples.

\noindent\textbf{Testing:}
All methods are evaluated on the same 50 test cases at $t=50$ that
are not used for training. Each test case corresponds to one
basal-friction input sample. For each ICL prediction, five coarse--fine
demonstrations are selected from the 20 paired training samples, and
the coarse solution of the test case is used as the query. The
model predicts the corresponding fine velocity field without
any parameter updates.

\paragraph{Results.}

We evaluate all methods on the same 50 test cases at $t=50$. All
three learned methods use the same 20 high-fidelity velocity training
samples. ICL and the non-composite model use 20 paired low-fidelity
training samples, while the composite model uses all 100
low-fidelity training samples. The results are summarized in
Table~\ref{tab:sec34_reference_mf}.

\begin{table}[H]
\centering
\small
\setlength{\tabcolsep}{7pt}
\renewcommand{\arraystretch}{1.15}

\begin{tabular}{lccccc}
\toprule

Method
& Test cases
& \multicolumn{3}{c}{
$\displaystyle
\frac{\|\mathbf{v}_c-\mathbf{v}_f\|_2}
{\|\mathbf{v}_f\|_2}$ (\%)
}
& MSE $(\mathrm{m}/\mathrm{yr})^2$ \\
\cmidrule(lr){3-5}

&
&
Mean
& Median
& P90
& \\

\midrule

LF baseline
& 50
& 27.35
& 27.60
& 29.84
& 14.302 \\

\midrule

Method
& Test cases
& \multicolumn{3}{c}{
$\displaystyle
\frac{\|\widehat{\mathbf{v}}_f-\mathbf{v}_f\|_2}
{\|\mathbf{v}_f\|_2}$ (\%)
}
& MSE $(\mathrm{m}/\mathrm{yr})^2$ \\
\cmidrule(lr){3-5}

&
&
Mean
& Median
& P90
& \\

\midrule

Non-composite MF
& 50
& 7.70
& 7.03
& 12.15
& 2.171 \\

Composite MF
& 50
& 7.14
& 6.84
& 10.95
& 1.787 \\

ICL
& 50
& \textbf{4.93}
& \textbf{3.68}
& \textbf{10.48}
& \textbf{1.262} \\

\bottomrule
\end{tabular}

\caption{
Multi-fidelity prediction on 50 ice-sheet test cases at $t=50$ that
are not used for training. All learned methods use the same 20
high-fidelity velocity training samples. ICL and the non-composite
model use 20 paired low-fidelity samples, while the composite model
uses 100 low-fidelity samples. The upper block reports the
coarse--fine discrepancy of the low-fidelity baseline, while the
lower block reports the prediction errors of the learned
multi-fidelity models. Both blocks report the mean, median, and
90th-percentile relative $L^2$ error together with the MSE.
}
\label{tab:sec34_reference_mf}
\end{table}

\FloatBarrier

\paragraph{Results analysis.}

The LF baseline has a mean relative $L^2$ error of $27.35\%$. The
non-composite and composite multi-fidelity DeepONets reduce the mean
error to $7.70\%$ and $7.14\%$, respectively, while ICL achieves the
lowest mean error of $4.93\%$. The median ICL error is $3.68\%$, the
90th-percentile error is $10.48\%$, and the mean coefficient of
determination is $R^2=0.9960$.

Compared with the non-composite multi-fidelity model, ICL reduces the
mean prediction error by approximately $36.0\%$ while using the same
20 low-fidelity and 20 high-fidelity velocity training samples.
Compared with the composite model, ICL reduces the mean prediction
error by approximately $30.9\%$ and the MSE from $1.787$ to $1.262$,
while using the same 20 high-fidelity velocity training samples and
only 20 low-fidelity training samples instead of 100. Thus, ICL
achieves the best prediction accuracy while using 80 fewer
low-fidelity training samples than the composite multi-fidelity model.


\subsubsection{Generalization to unseen operators}

We next extend the single-operator reference benchmark to a
multi-operator ICL setting and examine whether the model can
generalize to ice-sheet operators that are unseen during training.
Both the rheological parameter $A$ and the basal-friction field
$\beta(x,y)$ vary across operators, and each fixed operator provides
solution states at different times.

For each prediction, the model is provided with coarse--fine
demonstrations from one unseen testing operator at different times.
A new coarse solution $\mathbf v_{c,q}$ from another time of that
operator is used as the query. The model receives the spatial grid,
time, demonstrations, and coarse query, but not the operator-defining
quantities $A$ and $\beta(x,y)$. The model uses the demonstrations
as context to predict $\widehat{\mathbf v}_{f,q}$ for the unseen
operator.

\noindent\textbf{Training:}
We generate 70 training operators by varying $A$ and independently
sampling the basal-friction field $\beta(x,y)$. Each operator contains
paired coarse and fine velocity fields at 50 time points.

Each training example is constructed from a single operator. The
demonstrations and query are taken from different time points of that
operator, and the number of demonstrations is varied dynamically from
$D=1$ to $D=5$ during training.

\noindent\textbf{Testing:}
We use 30 additional operators as unseen testing operators. These
operators are generated from the same operator distribution as the
training operators and are not used during training. For each unseen
testing operator, five time points are used as coarse--fine
demonstrations, and each of the remaining 45 time points defines one
test case. This gives $30\times45=1,350$ test cases in total.

\paragraph{Results.}

We evaluate the model on 1,350 test cases from the 30 unseen testing
operators. Each prediction uses five coarse--fine demonstrations from
the same testing operator. We also evaluate a wrong-operator context
in which the query is unchanged but all five demonstrations are
replaced by coarse--fine pairs from another unseen testing operator.
The results are summarized in Table~\ref{tab:sec34_unseen_operators}.

\begin{table}[H]
\centering
\small
\setlength{\tabcolsep}{9pt}
\renewcommand{\arraystretch}{1.15}

\begin{tabular}{lcccc}
\toprule

Evaluation
& Test cases
& \multicolumn{3}{c}{
$\displaystyle
\frac{\|\mathbf{v}_c-\mathbf{v}_f\|_2}
{\|\mathbf{v}_f\|_2}$ (\%)
}
\\
\cmidrule(lr){3-5}

&
&
Mean
& Median
& P90
\\

\midrule

LF baseline
& 1350
& 30.62
& 30.85
& 33.27
\\

\midrule

Evaluation
& Test cases
& \multicolumn{3}{c}{
$\displaystyle
\frac{\|\widehat{\mathbf{v}}_f-\mathbf{v}_f\|_2}
{\|\mathbf{v}_f\|_2}$ (\%)
}
\\
\cmidrule(lr){3-5}

&
&
Mean
& Median
& P90
\\

\midrule

ICL, correct context
& 1350
& \textbf{1.15}
& \textbf{0.95}
& \textbf{1.86}
\\

ICL, wrong-operator context
& 1350
& 1.49
& 1.06
& 2.63
\\

\bottomrule
\end{tabular}

\caption{
Unseen-operator prediction on 1,350 test cases from 30 unseen testing
ice-sheet operators using five demonstrations. The upper block reports
the coarse--fine discrepancy of the low-fidelity baseline, while the
lower block reports the ICL prediction error with the correct and
wrong-operator contexts. In the wrong-operator context, the query is
unchanged while all five demonstrations are replaced by coarse--fine
pairs from another unseen testing operator. Both blocks report the
mean, median, and 90th-percentile relative $L^2$ error.
}
\label{tab:sec34_unseen_operators}
\end{table}

\begin{figure}[H]
    \centering
    \includegraphics[width=\textwidth]{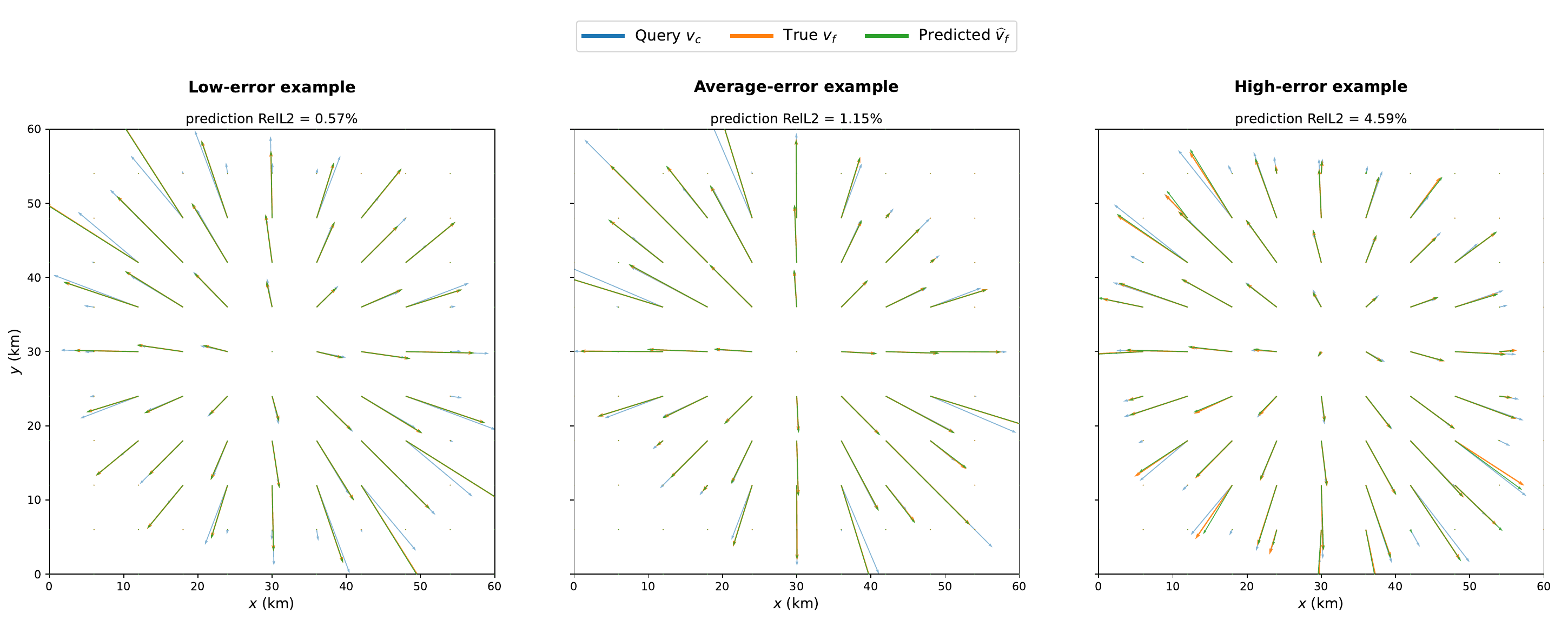}
    \caption{
Representative predictions for unseen testing two-dimensional
ice-sheet operators. The panels show low-, average-, and high-error
examples. The coarse query, fine target, and ICL prediction are
overlaid as blue, orange, and green velocity vectors, respectively.
The relative $L^2$ prediction error is reported above each panel.
    }
    \label{fig:sec34_ice_vars}
\end{figure}

Representative low-, average-, and high-error predictions are shown in
Figure~\ref{fig:sec34_ice_vars}. In each example, the low-fidelity
query velocity field, the corresponding fine target, and the
ICL prediction are overlaid to compare both the direction and
magnitude of the two-dimensional velocity vectors. The ICL relative
$L^2$ error reported in each panel measures
the discrepancy between the predicted and true high-fidelity velocity
fields.

\FloatBarrier

\paragraph{Results analysis.}

Across the 1,350 test cases, the LF baseline has a mean relative
$L^2$ error of $30.62\%$, while ICL reduces the mean prediction
error to $1.15\%$, corresponding to an approximately $96.2\%$
reduction. The median prediction error is $0.95\%$, and the
90th-percentile error is $1.86\%$. The mean coefficient of
determination is $R^2=0.9998$.

These results show that a single shared ICL model can accurately
predict the fine velocity fields for operators that are not
observed during training.

The prediction also depends on the operator-specific
demonstrations. Replacing the five correct demonstrations with
demonstrations from another unseen testing operator increases the mean
error from $1.15\%$ to $1.49\%$, corresponding to a relative increase
of approximately $30.1\%$. The median error increases from $0.95\%$
to $1.06\%$, while the 90th-percentile error increases from $1.86\%$
to $2.63\%$.

The model is not given the operator-defining quantities $A$ or
$\beta(x,y)$. Instead, it receives the spatial grid, time, the
coarse--fine demonstrations, and the coarse query. The increase in
error under the wrong-operator context shows that the model uses
operator-specific information supplied by the demonstrations while
adapting to previously unseen ice-sheet operators without parameter
updates.

\section{Discussion and conclusion}

In this work, we address the multi-fidelity problem through in-context learning. Multi-fidelity learning combines inexpensive low-fidelity solutions with limited high-fidelity data to improve prediction accuracy while reducing computational cost. However, many multi-fidelity operator-learning methods are developed for a specific coarse-to-fine relation or a single operator-learning task, making extension to multiple operators costly. We propose an in-context multi-fidelity learning framework that learns a family of coarse-to-fine operators with one shared model using a few coarse--fine demonstrations to provide the information needed for the trained model to predict the fine solution for a new coarse query. Once the model is trained, its parameters remain fixed, and no fine-tuning is required at inference.

The numerical experiments show that the framework generalizes to unseen
operators across several PDE problems and remains effective for
operator-parameter and query-side OOD prediction. The model uses the
supplied demonstrations, with prediction accuracy depending on both
their number and consistency. The framework also handles temporally
shifted coarse-to-fine mappings and extends to a two-dimensional
multi-fidelity problem.

The cross-operator experiments show that OOD information learned from
operators with more training data can improve prediction for operators
with limited data that have not observed the same OOD regime. Previously available demonstrations can be reused for these OOD queries.
As a result, a new OOD prediction requires only a new coarse query
generated by the low-fidelity solver and does not require new
high-fidelity data or fine-tuning. This is particularly advantageous in
multi-fidelity settings, where low-fidelity simulations are inexpensive
while high-fidelity simulations are costly.


\section*{Acknowledgments}
ZZ and XL are supported by DOE award number DE-SC0025440 and NSF award number 2615564.
The work of PS was completed with support from the U.S. Department of Energy's Advanced Scientific Computing Research program under the Uncertainty Quantification for Multifidelity Operator Learning (MOLUcQ) project (Project No. 81739). Pacific Northwest National Laboratory (PNNL) is a multi-program national laboratory operated for the U.S. Department of Energy (DOE) by Battelle Memorial Institute under Contract No. DE-AC05-76RL01830. MZ was supported by AFOSR award number FA9550-22-1-0004 and NSF award number CBET-2338843.


\appendix

\section{Additional experimental details}
\label{app:experimental_details}

\subsection{Common implementation details}
\label{app:common_details}

\subsubsection{Variable-length Transformer details}
\label{app:transformer_details}

For a valid input sequence
$S_i^D=[r_1,\ldots,r_{L_D}]^\top\in\mathbb R^{L_D\times N_r}$,
let $\mathbf h_\ell\in\mathbb R^{d_{\mathrm{model}}}$ denote
the embedding of row $r_\ell$ supplied to the Transformer,
where $d_{\mathrm{model}}$ is the model dimension. The embeddings
are constructed using the application-specific input encoders
described in Appendices~\ref{app:scalar_details}--\ref{app:ice_details}. They form
$H_i^D=[\mathbf h_1,\ldots,\mathbf h_{L_D}]^\top
\in\mathbb R^{L_D\times d_{\mathrm{model}}}$.
For one
self-attention head, let
$\mathbf W_Q,\mathbf W_K,\mathbf W_V
\in\mathbb R^{d_{\mathrm{model}}\times d_k}$ denote the query, key,
and value projection matrices. Then
$\mathbf Q_i^D=H_i^D\mathbf W_Q$,
$\mathbf K_i^D=H_i^D\mathbf W_K$, and
$\mathbf V_i^D=H_i^D\mathbf W_V$, with
$\mathbf Q_i^D,\mathbf K_i^D,\mathbf V_i^D
\in\mathbb R^{L_D\times d_k}$ and
$\mathbf Q_i^D(\mathbf K_i^D)^\top
\in\mathbb R^{L_D\times L_D}$. Changing the demonstration count $D$
changes the sequence and attention dimensions through $L_D$, but does
not change the dimensions of the trainable projection matrices.

\subsubsection{Training objective}
\label{app:training_objective}

For a training batch of size $B$, let
$\widehat{\mathbf y}_b,\mathbf y_b\in\mathbb R^{N_r}$
denote the predicted and target fine query rows for training example
$b=1,\ldots,B$. The loss is evaluated using the normalized rows,
with application-specific normalization described in Appendices~\ref{app:scalar_details}--\ref{app:ice_details}.
The training objective is the mean-squared error
\[
\mathcal L(\theta)
=
\frac{1}{B N_r}
\sum_{b=1}^{B}
\sum_{\ell=1}^{N_r}
\left(\widehat y_{b,\ell}-y_{b,\ell}\right)^2.
\]

\subsubsection{Batching strategy}
\label{app:batching_strategy}

For experiments with a dynamic number of demonstrations, the
demonstration count $D$ is sampled independently for each training
example. Examples with different values of $D$ can have different
valid sequence lengths $L_D$. Within each training batch, padded rows are appended after the final
valid row of shorter sequences until the application-specific maximum
sequence length is reached, and an attention padding mask prevents the
padded rows from contributing to the valid attention context.
Padding rows are excluded when determining $L_D$.
The output at the final valid input position $L_D$, corresponding
to the coarse-query row, is used to predict the missing next row
$\widehat r_{L_D+1}=\widehat{\mathbf v}_{f,q}^{(i)}$.
The maximum padded sequence length used for each application is
reported in Table~\ref{tab:app_sequence_summary}.

\subsubsection{Input-sequence summary}
\label{app:sequence_summary}

Across all applications, demonstration $j$ follows the row structure
$[\boldsymbol{\alpha}_{d_j}^{(i)},
\mathbf v_{c,d_j}^{(i)},\mathbf v_{f,d_j}^{(i)}]$, while the query
follows
$[\boldsymbol{\alpha}_q^{(i)},\mathbf v_{c,q}^{(i)}]$.
As defined in Section~\ref{sec:icl_framework}, $a(u_s^{(i)})$ is the application-specific
information for sample $s$, and $\boldsymbol{\alpha}_s^{(i)}$
expands into its $n_a$ auxiliary row vectors
$\mathbf a_{s,1}^{(i)},\ldots,\mathbf a_{s,n_a}^{(i)}$,
each in $\mathbb R^{N_r}$ after any required interpolation or
broadcasting. Table~\ref{tab:app_sequence_summary} lists the
quantities represented by these rows in sequence order.
The valid input length is $L_D=D(n_a+2)+(n_a+1)$, and the
maximum padded length corresponds to $D=5$.

\begin{table}[H]
\centering
\small
\setlength{\tabcolsep}{5pt}
\renewcommand{\arraystretch}{1.15}
\begin{tabular}{lcccc}
\toprule
Application & Auxiliary-row quantities & $n_a$ & $L_D$ & Max. length \\
\midrule
Scalar conservation laws & $[x]$ & 1 & $3D+2$ & 17 \\
KS & $[x,t]$ & 2 & $4D+3$ & 23 \\
KdV & $[x]$ & 1 & $3D+2$ & 17 \\
Ice sheet, reference & $[XY,\beta,H]$ & 3 & $5D+4$ & 29 \\
Ice sheet, multi-operator & $[XY,t]$ & 2 & $4D+3$ & 23 \\
\bottomrule
\end{tabular}
\caption{Application-specific input-sequence structure. The maximum
length corresponds to $D=5$ demonstrations. For KdV, the demonstrated
fine row is evaluated at the shifted target time $t+\Delta t$.}
\label{tab:app_sequence_summary}
\end{table}


\subsection{Scalar conservation-law experiments}
\label{app:scalar_details}

\paragraph{Data construction and numerical solvers.}

For the standard unseen-operator experiments, each of the eight PDE
families contains 100 operators. For each operator, 100 solution
samples are generated by varying the initial condition. The operator
parameters are sampled as $q_1\sim\mathcal U(0.5,2.0)$ for all
families and $q_2\sim\mathcal U(0.005,0.02)$ for the viscous
families, while $q_2=0$ for the inviscid families.

The initial condition is a smooth periodic box profile,
\[
u_0(x)
=
\frac{A}{2}
\left[
\tanh\left(\frac{\xi+w/2}{\epsilon}\right)
-
\tanh\left(\frac{\xi-w/2}{\epsilon}\right)
\right],
\qquad
\xi=((x-c+0.5)\bmod1)-0.5,
\]
where $c\sim\mathcal U(0,1)$, $w\sim\mathcal U(0.10,0.30)$,
$A\sim\mathcal U(0.80,1.50)$, and
$\epsilon\sim\mathcal U(0.015,0.035)$.

All coarse and fine solutions used by the model are evaluated at the
terminal time $t=2$. The coarse solution is computed on a
50-point grid using the low-fidelity HLL-based finite-volume solver with shock
viscosity. The fine solution is computed on a 100-point grid
using the high-fidelity MUSCL--Hancock solver with a Superbee limiter. The
coarse solution is linearly interpolated from 50 to 100 points
before being provided to the model, so that the coarse and fine
solutions use the same spatial representation.

\paragraph{Unseen-operator training and evaluation.}

Within each PDE family, 70 operators are used for training and 30 are
used as unseen testing operators, giving 560 training operators and
240 unseen testing operators in total. Normalization statistics are
computed using the training operators.

Each training example is constructed from a single operator. The
demonstrations and query are generated from different initial
conditions of that same operator.

The row-sequence Transformer has model dimension 96, four attention
heads, three encoder layers, feed-forward dimension 192, and
dropout 0.1, giving a total of 255,364 trainable parameters. Training
uses 300 epochs with 1,200 optimization steps per epoch and batch size
1,000. AdamW is used with weight decay $10^{-4}$. The learning rate
increases from $10^{-4}$ to $2\times10^{-4}$ during an 8\% warmup
period and then follows a cosine schedule to $10^{-5}$.

The unseen-operator evaluation uses 1,024 fixed test cases.
Each test case contains five coarse--fine demonstrations from one
unseen testing operator and a coarse query generated from another
initial condition of the same operator. The trained model directly
predicts the corresponding fine solution for each query.

\paragraph{Scaling and context interventions.}

For the demonstration-scaling experiment, the same fixed queries are
evaluated using nested contexts with $D=1,\ldots,5$ demonstrations. The same trained model is used for all values of $D$.

For the wrong-operator intervention, all five demonstrations are
replaced by coherent coarse--fine pairs from another unseen testing operator
in the same PDE family, while the original query is retained. For the
wrong-$\mathbf v_f$ intervention, the original demonstrated coarse solutions
are retained while their fine counterparts are replaced by
those from another unseen testing operator in the same PDE family.
These interventions use the same
fixed test queries as the standard five-demonstration evaluation.

\paragraph{Operator-parameter OOD evaluation.}

The operator-parameter OOD experiment evaluates the same trained model
on new operators outside the training range of $q_1$.

For each PDE family, 50 new operators are generated with
$q_1\in[0.45,0.50)$ and another 50 with
$q_1\in(2.00,2.20]$, giving 800 OOD operators in total. The $q_2$
setting, initial-condition distribution, and numerical solvers remain
the same as in the standard experiment.

For each OOD test case, five coarse--fine demonstrations and one
query are generated from the same OOD operator using different initial
conditions. The evaluation contains 1,024 fixed test cases, with
498 cases from the lower-OOD region and 526 from the upper-OOD region.
The same trained model directly predicts the corresponding fine
solution for each query.

\paragraph{Initial-condition amplitude OOD experiment.}

The amplitude-OOD experiment uses a separate dataset containing 200
operators for each of the eight PDE families, giving 1,600 operators
in total. The operator parameters $q_1$ and $q_2$ are generated from
the same ranges used in the standard unseen-operator experiment.
Within each family, 100 operators are assigned to Group A and 100 to
Group B, resulting in 800 operators in each group.

The base-amplitude range is $A\in[0.80,1.50]$, while the near- and
farther-OOD ranges are $A\in(1.50,1.57]$ and $A\in(1.57,1.65]$,
respectively.

Each operator contains 100 matched initial-condition samples. For a
matched sample across the three amplitude ranges, the center $c$,
width $w$, and edge-smoothing parameter $\epsilon$ are kept fixed,
while the amplitude $A$ is changed.

We train two models using the same architecture, normalization
procedure, optimizer settings, and training budget. The models differ
in the amplitude ranges included in their training data. For both
models, normalization statistics are computed using the base-amplitude
data from all 1,600 operators.

For the \emph{base-amplitude-only} model, both Group A and Group B are
trained on samples from $A\in[0.80,1.50]$.

For the \emph{mixed-amplitude} model, Group A is trained on the
base-amplitude range $A\in[0.80,1.50]$, the near-OOD range
$A\in(1.50,1.57]$, and the farther-OOD range $A\in(1.57,1.65]$.
Group B is trained on the base-amplitude range
$A\in[0.80,1.50]$.

Both models are evaluated on the same 800 Group B operators. For each
operator, all 100 matched initial-condition samples are used to
construct near- and farther-OOD queries, giving
$800\times100=80,000$ test cases in each OOD regime.

For each Group B query, five base-amplitude coarse--fine pairs from
the same Group B operator are reused as demonstrations. The query
sample itself is excluded from the demonstration set, and the
base-amplitude-only and mixed-amplitude models use exactly the same
five demonstrations for each matched query. The selected
demonstrations are base-amplitude pairs included in training. A new
coarse query $\mathbf v_{c,q}$ is generated from the
corresponding matched initial condition with its amplitude changed to
either the near- or farther-OOD range. The trained model then predicts
the corresponding fine solution.

Group B predictions use demonstrations from the same Group B
operator. In the mixed-amplitude setting, the OOD amplitude
information introduced through Group A is learned by the shared model
and can improve prediction for Group B. The corresponding Group B
high-fidelity OOD solution is used only as the reference solution for
computing the evaluation error.


\subsection{Kuramoto--Sivashinsky equation}
\label{app:ks_details}

\paragraph{Data construction and numerical solvers.}

For the unseen-operator KS experiment, we generate 200 operators with
$\beta$ distributed over $[0.070,0.085]$. For each operator, 20
independent solution trajectories are generated over $t\in[0,3.5]$.
The numerical solver time step is $0.005$, and solution snapshots are
stored every $0.05$ time units, giving 71 saved time points per
trajectory.

Each trajectory uses the initial condition
$u_0(x)=f\!\left(\sin^2(x/2)\right)$, where $f$ is a zero-mean
squared-exponential Gaussian random field with covariance
$\operatorname{Cov}[f(z),f(z')]=\sigma^2\exp\!\left(-(z-z')^2/(2\ell^2)\right)$,
with $\ell=0.5$ and $\sigma=1.0$.
The 20 trajectories for each operator use independent realizations of
this random field. The initial state is first constructed on the
256-point high-fidelity grid, and the same underlying numerical state
is periodically Fourier-resampled to the 72-point low-fidelity grid
and the 100-point common learning grid.

The coarse solution is computed on a 72-point periodic grid using
the low-fidelity second-order centered finite-difference spatial discretization with
ETDRK4 time integration. The fine solution is computed on a
256-point periodic grid using the high-fidelity Fourier pseudospectral discretization
with $3/2$ dealiasing and ETDRK4 time integration.

Before being provided to the model, the low- and high-fidelity
solutions are represented on a common 100-point periodic learning grid.
The resulting low- and high-fidelity arrays have dimensions
$200\times20\times71\times100$, corresponding to operator, trajectory,
saved time, and spatial location.

Although each trajectory is simulated over $t\in[0,3.5]$, only the
late-time interval $t\in[2.50,3.50]$ is used to construct the ICL
training examples. With snapshots separated by $0.05$, this interval contains
21 eligible time points.

\paragraph{Unseen-operator training and evaluation.}

The 200 operators are divided using a fixed operator-level split into
140 training operators and 60 unseen testing operators. Normalization
statistics are computed using only the training operators and the
eligible time points in $t\in[2.50,3.50]$.

Each training example is constructed from a single training operator
and a single trajectory of that operator. The demonstrations and query
are selected from distinct eligible time points of the same trajectory.
The prediction target is $\mathbf v_f(t_q)$ at the query time, and the
operator parameter $\beta$ is not provided explicitly to the model.

The KS model uses a row-sequence Transformer with model dimension 288,
six attention heads, five encoder layers, feed-forward dimension 576,
and dropout 0.1, giving 3,484,036 trainable parameters. Training uses
a batch size of 1,000 for 150 epochs with 900 optimization steps per epoch.
AdamW is used with weight decay $10^{-4}$. The learning rate increases
from $10^{-4}$ to $2\times10^{-4}$ during an 8\% warmup period and
then follows a cosine schedule to $10^{-5}$. The best checkpoint is
selected using a fixed 512-case $D=5$ monitor constructed only from the
training operators.

The final unseen-operator evaluation uses 960 fixed test cases from the
60 unseen testing operators, with 16 cases constructed for each
operator. Each test case contains five coarse--fine demonstrations and
one distinct query selected from the same trajectory of one testing
operator, with all demonstration and query times drawn from
$t\in[2.50,3.50]$. The trained model directly predicts the corresponding
fine solution at the query time without any parameter updates.

\paragraph{Scaling and context interventions.}

For the demonstration-scaling evaluation, the same 960 fixed test cases
are evaluated using nested contexts with $D=1,\ldots,5$
demonstrations. The context with $D$ demonstrations contains the first
$D$ demonstrations of the corresponding five-demonstration test case.
Thus, the query is held fixed while only the amount of valid
demonstration context is changed. The same trained model is used for
all values of $D$, with no parameter updates.

For the context-intervention evaluations, the same 960 fixed test cases
and five-demonstration setting are used. In the correct-context
condition, the original five coarse--fine demonstration pairs from the
same operator and trajectory as the query are retained. In the
wrong-operator condition, all five demonstration pairs are replaced by
coherent coarse--fine pairs from another unseen operator at the
corresponding demonstration times, while the original query is left
unchanged. In the wrong-$\mathbf v_f$ condition, the original
demonstrated coarse states are retained, but the corresponding
fine demonstration states are replaced by fine states
from the donor operator. This preserves the original query and coarse
demonstration inputs while disrupting the demonstrated coarse-to-fine
correspondence.

The mean prediction errors for $D=1,\ldots,5$ are $9.37\%$, $6.63\%$,
$5.51\%$, $4.85\%$, and $4.39\%$, respectively. For the
context-intervention study, the mean prediction errors are $4.39\%$
with the correct context, $52.10\%$ with the wrong-operator context,
and $127.28\%$ with the wrong-$\mathbf v_f$ context.

\paragraph{Operator-parameter OOD evaluation.}

The operator-parameter OOD experiment uses the same trained model,
normalization statistics, numerical solvers, and late-time evaluation
interval as the standard unseen-operator KS experiment. No OOD operator
is used for training, normalization, training-side monitoring, or
checkpoint selection.

We generate 50 unseen OOD operators with
$\beta\in(0.085,0.100]$. Each operator contains 20 trajectories
constructed using the same initial-condition procedure, numerical
solvers, common learning grid, and saved time points as in the standard
experiment.

The final evaluation contains 512 fixed test cases. Each test case uses
five same-operator, same-trajectory coarse--fine demonstrations and one
distinct query with times in $t\in[2.50,3.50]$. The corresponding fine
query state is used only as the evaluation target.

\paragraph{Cross-operator temporal OOD experiment.}

The cross-operator temporal OOD experiments use separate KS datasets
with the same 200 uniformly spaced $\beta$ values in
$[0.070,0.085]$ and 20 independent trajectories for each operator. The
same low- and high-fidelity spatial discretizations are used as in the
unseen-operator experiment: the low-fidelity solver uses the 72-point
second-order centered finite-difference discretization with ETDRK4, the
high-fidelity solver uses the 256-point $3/2$-dealiased Fourier
pseudospectral discretization with ETDRK4, and both fidelities are
represented on the common 100-point learning grid.

Both experiments use the same base-time pool of 100 points,
$2.505,\,2.510,\,\ldots,\,3.000$. For the evaluation interval
$[3.002,3.100]$, the temporal OOD pool contains the 50 points
$3.002,\,3.004,\,\ldots,\,3.100$, separated by $0.002$.
For the evaluation interval $[3.005,3.250]$, the temporal OOD pool
contains the 50 points $3.005,\,3.010,\,\ldots,\,3.250$,
separated by $0.005$; the numerical solver time step for this dataset
is $0.001$. Each experiment has low- and high-fidelity arrays of
dimensions $200\times20\times150\times100$.

A fixed random split with seed 42 divides the 200 operators into
100 Group A operators and 100 Group B operators. For each evaluation
interval, we train two matched models using the same architecture,
optimization settings, normalization procedure, and training budget.
The two models differ only in the temporal data available to Group A.

For the \emph{base-time-only} model, both Group A and Group B are
trained using only the 100 base-time points. For the
\emph{mixed-time} model, Group A is trained using its 100 base-time
points together with the 50 temporal OOD points, while Group B is
trained using only the 100 base-time points. Thus, the Group A
mixed-time eligible pool contains 150 time points, of which one third
are from the temporal OOD interval. Group B OOD states are never used
for training or checkpoint selection.

Each training example is constructed from one fixed $\beta$ operator
and one trajectory of that operator, with the demonstrations and query
selected from distinct eligible time points.

Both cross-operator models use the same KS Transformer architecture as
the standard experiment: model dimension 288, six attention heads,
five encoder layers, feed-forward dimension 576, dropout 0.1, and
3,484,036 trainable parameters. Training uses a batch size of 1,000,
150 epochs, and 900 optimization steps per epoch with AdamW, weight
decay $10^{-4}$, an 8\% warmup from $10^{-4}$ to $2\times10^{-4}$,
and cosine decay to $10^{-5}$. Normalization is matched between the two
models and is based on the shared base-time data. Checkpoint selection
uses matched base-time training-side monitoring and does not use Group B
OOD data.

The final paired evaluation covers all 100 Group B operators, all 20
trajectories per operator, and all 50 temporal OOD query times, giving
$100\times20\times50=100{,}000$ test cases for each evaluation interval. Every test case uses
$D=5$ base-time demonstrations from the same Group B operator and the
same trajectory as the query. The evaluation uses demonstration
contexts that are shared between the two trained models, so the
base-time-only and mixed-time predictions are compared on exactly the
same queries and demonstrations. The corresponding Group B
high-fidelity OOD states are used only as evaluation targets.

For query times in $[3.002,3.100]$, the mean coarse--fine discrepancy
is $28.77\%$. The base-time-only and mixed-time models have mean
prediction errors of $2.90\%$ and $2.70\%$, respectively. The
relative improvement is $6.72\%$, and the paired win rate is
$57.71\%$. The corresponding mean $R^2$ values are $0.9982$ and
$0.9984$, respectively.

For query times in $[3.005,3.250]$, the mean coarse--fine discrepancy
is $30.52\%$. The base-time-only and mixed-time models have mean
prediction errors of $6.62\%$ and $5.08\%$, respectively. The
relative improvement is $23.33\%$, and the paired win rate is
$69.62\%$. The corresponding mean $R^2$ values are $0.9865$ and
$0.9919$, respectively.


\subsection{Korteweg--de Vries equation}
\label{app:kdv_details}

\paragraph{Data construction and numerical solvers.}

The initial condition is $u(x,0)=A\cos(x+\phi)$, where
$A\sim\mathcal U(0.35,0.65)$ and
$\phi\sim\mathcal U(0,2\pi)$. The dispersion parameter is sampled as
$\epsilon\sim\mathcal U(0.01,0.02)$. One operator is defined by the
fixed triple $(\epsilon,A,\phi)$.

We generate 200 independent operators. Each operator is simulated over
$t\in[0,10]$. The coarse solution is computed on a 32-point
periodic spatial grid, while the fine solution is computed on
a 256-point periodic spatial grid. The trajectories are generated
using Fourier spectral discretization with ETDRK4 time integration and
a solver time step of $0.002$. Solution snapshots are stored every
$0.05$ time units.

Before being provided to the model, the low- and high-fidelity
solutions are independently represented on a common 100-point periodic
grid using piecewise-linear interpolation.

For each operator, 181 paired source--target samples are constructed. The source times span $t\in[0,9]$, while each corresponding fine target is evaluated at $t+0.5$, so that every paired sample has the form $(\mathbf v_c(t),\mathbf v_f(t+0.5))$.

\paragraph{Unseen-operator training and evaluation.}

The 200 operators are divided using a fixed operator-level split into
140 training operators and 60 unseen testing operators. Normalization
statistics are computed using only the 140 training operators.

Each training example is constructed from a single training operator.
The demonstrations and query are selected from distinct source-time
samples of that same operator.

The Transformer has model dimension 256, four attention heads, five
encoder layers, feed-forward dimension 512, and dropout 0.1, giving
2,759,012 trainable parameters. Training uses a batch size of 1,000
for 150 epochs with 1,200 optimization steps per epoch. AdamW is used
with weight decay $10^{-4}$ and gradient clipping at 1.0. The learning
rate increases from $10^{-4}$ to $2\times10^{-4}$ during an $8\%$
warmup period and then follows a cosine schedule to $10^{-5}$.

The final unseen-operator evaluation uses 960 fixed test cases from
the 60 unseen testing operators, with 16 cases per operator. Each test
case contains five distinct demonstrations and one distinct query from
the same testing operator.

\paragraph{Scaling and context interventions.}

For the demonstration-scaling evaluation, the same 960 fixed queries
are evaluated using nested contexts with $D=1,\ldots,5$
demonstrations. The context with $D$ demonstrations contains the first
$D$ demonstrations of the corresponding five-demonstration test case.
The same trained model is used for all values of $D$.

For the wrong-operator intervention, all five coarse--fine
demonstration pairs are replaced by pairs from another unseen testing
operator at the same demonstration source times, while the original
query is retained. For the wrong-$\mathbf v_f$ intervention, the
original demonstrated coarse states are retained while their
fine counterparts are replaced by the corresponding states
from another unseen testing operator. Both interventions use the same
960 fixed queries as the standard five-demonstration evaluation.


\subsection{Ice-sheet flow experiments}
\label{app:ice_details}

\paragraph{Reference benchmark and data construction.}

The physical parameters, constitutive settings, Halfar-dome initial
condition, and basal-friction construction are motivated by the
reference multi-fidelity benchmark \cite{howard2023multifidelity},
with the modifications described below. All
numerical solutions used in this study are regenerated independently
using our reduced ice-flow solver rather than taken from the original
benchmark data.

The basal-friction perturbation is generated as
$\log\beta=\log\beta_{\mathrm{opt}}+0.35\,\eta$. The reference
basal-friction field and spatial correlation construction follow
the reference benchmark \cite{howard2023multifidelity}.

The coarse and fine velocity solutions use the same reduced
governing equations at different spatial resolutions. The
low-fidelity grid contains $15\times15$ points, while the
high-fidelity grid contains $41\times41$ points. The low-fidelity
velocity field is bilinearly interpolated to the common
$41\times41$ representation before being provided to the model.

For the reference experiment, the coarse and fine thickness and
velocity trajectories are generated separately on the $15\times15$
and $41\times41$ grids, respectively, using the same underlying
basal-friction realization for each paired sample.
The low-fidelity thickness evolves on its own solver grid rather
than being obtained by downsampling the high-fidelity thickness.
For the remaining 80 low-fidelity training samples without paired
high-fidelity data, only low-fidelity trajectories are generated.

The ice thickness $H$ is part of the application-specific information
$a(u)$ associated with the velocity solution. The ICL model uses
the high-fidelity-grid basal-friction field $\beta$ and ice thickness
$H$, together with the interpolated coarse velocity.
Both the composite and non-composite multi-fidelity DeepONet
baselines also use the high-fidelity-grid $\beta$ and $H$
for high-fidelity prediction.

\paragraph{Numerical solver.}

Both fidelity levels use the same reduced membrane-style governing
equations, with fidelity determined by the spatial grid resolution.
The nonlinear velocity problem is solved using Picard iteration with
a relative convergence tolerance of $10^{-5}$ and at most 30
iterations per block. Adaptive under-relaxation uses factors
$0.65$, $0.50$, $0.35$, and $0.25$ as needed for convergence.

The ice-thickness equation is advanced using a conservative upwind
scheme with a CFL number of $0.45$ and a maximum internal time step of
$0.20$ years. For the trajectory-based multi-operator dataset, each
simulation is advanced from $t=0$ to $t=50$, and solution fields are
saved at $t=1,2,\ldots,50$.

\paragraph{Reference multi-fidelity benchmark.}

The reference experiment uses a single-operator multi-fidelity
setting motivated by the reference benchmark
\cite{howard2023multifidelity}. We use 100
low-fidelity training samples and 20 high-fidelity velocity training
samples at the final time $t=50$. The 20 high-fidelity velocity
samples are paired with their corresponding coarse solutions,
giving 20 paired coarse--fine training samples. The test set contains
50 test cases available at both fidelity levels, with each test case
corresponding to one basal-friction input sample at $t=50$.

The ICL training set consists of the 20 paired coarse--fine samples.
Demonstrations and queries are formed from different paired samples
under the same reference setting. At test time, five paired training
samples are used as demonstrations for each test query.

The ICL model uses separate encoders for $XY$, $\beta$, and $H$,
together with a shared velocity encoder for $\mathbf v_c$ and
$\mathbf v_f$. The Transformer has model dimension 512, eight
attention heads, five encoder layers, feed-forward dimension 1536,
and dropout 0.1, with separate output heads for the two velocity
components.

Training uses a batch size of 64 for 180 epochs with 600 optimization
steps per epoch. AdamW is used with weight decay $2\times10^{-4}$.
The learning rate increases from $1.5\times10^{-5}$ to
$5\times10^{-5}$ during an 8\% warmup period and then follows a
cosine decay to $10^{-6}$. Normalization statistics are computed from
the 20 paired training samples.

The composite and non-composite multi-fidelity DeepONet baselines are
retrained on the same regenerated dataset following the corresponding
reference formulations \cite{howard2023multifidelity}. The composite model uses
100 low-fidelity training samples and 20 high-fidelity velocity
training samples. The non-composite model uses the 20 low-fidelity
samples paired with the 20 high-fidelity velocity training samples.
ICL uses the same 20 paired low- and high-fidelity velocity samples
as the non-composite model.

\paragraph{Multi-operator dataset construction.}

For the unseen-operator experiment, we extend the reference setting to
100 operators. Operator $i$ is defined by the fixed pair
$(A_i,\beta_i(x,y))$, where $A_i=A_0A_{\mathrm{rel},i}$ with
$A_0=10^{-16}~\mathrm{Pa^{-3}\,yr^{-1}}$. The 100 values of
$A_{\mathrm{rel},i}$ are evenly spaced over $[0.5,2.0]$.

Each operator is assigned an independent basal-friction field
$\beta_i(x,y)$ sampled from the modified basal-friction distribution
described in Appendix~\ref{app:ice_details}. The values of $A_i$ and $\beta_i(x,y)$ remain fixed
throughout the trajectory. The same $(A_i,\beta_i)$ pair is used for
the corresponding low- and high-fidelity simulations.

Each operator is simulated from $t=0$ to $t=50$, with paired coarse and
fine velocity fields saved at $t=1,2,\ldots,50$. This gives
50 paired time points per operator and 5,000 paired solution states in
total.

A fixed random operator-level split assigns 70 operators to training
and 30 operators as unseen testing operators. The training and unseen
testing operators are drawn from the same $A_{\mathrm{rel}}$ range
and basal-friction distribution. Normalization statistics and model
optimization are based on the 70 training operators.

\paragraph{Multi-operator ICL training.}

Each training example is constructed from a single training operator.
The demonstrations and query are selected from different time points
of that operator, with times drawn from the 50 available snapshots.

The multi-operator model uses separate encoders for $XY$ and scalar
time together with one shared velocity encoder for $\mathbf v_c$ and
$\mathbf v_f$. The Transformer has model dimension 128, four attention
heads, three encoder layers, feed-forward dimension 384, and dropout
0.1, with separate output heads for the two velocity components. This
gives 1,828,386 trainable parameters. Training uses the same
optimization schedule as the reference ICL model: batch size 64,
180 epochs, and 600 optimization steps per epoch with AdamW, an
8\% warmup, and cosine learning-rate decay. Its input representation
uses time to distinguish solution states along each operator trajectory.

\paragraph{Unseen-operator evaluation.}

For each of the 30 unseen testing operators, five fixed time points are used as
coarse--fine demonstrations. The remaining 45 time points are used as
queries, giving $30\times45=1,350$ test cases. For each query, the
five demonstrations come from the same unseen testing operator, and the trained
model predicts the corresponding fine velocity field.

For the wrong-operator intervention, the five demonstrations are
replaced by coarse--fine pairs from another unseen testing operator while the
query and its fine target remain fixed.


\bibliographystyle{unsrt}
\bibliography{references}

\end{document}